\documentclass[onefignum,onetabnum]{siamonline220329}

\usepackage{lipsum}
\usepackage{amsfonts}
\usepackage{graphicx}
\usepackage{epstopdf}
\usepackage{algorithmic}
\ifpdf
  \DeclareGraphicsExtensions{.eps,.pdf,.png,.jpg}
\else
  \DeclareGraphicsExtensions{.eps}
\fi

\usepackage{enumitem}
\setlist[enumerate]{leftmargin=.5in}
\setlist[itemize]{leftmargin=.5in}

\newsiamremark{remark}{Remark}
\newsiamremark{hypothesis}{Hypothesis}
\crefname{hypothesis}{Hypothesis}{Hypotheses}
\newsiamthm{claim}{Claim}

\headers{Stable adjusted root mean square propagation}{R. Li, J. Xu, and W. Xing}

\title{Stable gradient-adjusted root mean square propagation on least squares problem\thanks{Submitted to the editors DATE.
\funding{This work was funded in part by the National Natural Science Foundation of China Grant No. 11771243, PolyU postdoc matching fund scheme of the Hong Kong Polytechnic University Grant No. 1-W35A, and Huawei's Collaborative Grants ``Large scale linear programming solver'' and ``Solving large scale linear programming models for production planning''.}}}

\author{Runze Li \thanks{Department of Mathematical Sciences, Tsinghua University, Beijing 100084, China (\email{li-rz22@mails.tsinghua.edu.cn}).}
\and Jintao Xu\thanks{
Corresponding author. Department of Applied Mathematics, The Hong Kong Polytechnic University, Hong Kong, China (\email{xujtmath@163.com}).}
\and Wenxun Xing\thanks{Department of Mathematical Sciences, Tsinghua University, Beijing 100084, China (\email{wxing@tsinghua.edu.cn}).}}

\usepackage{amsopn}

\ifpdf
\hypersetup{
  pdftitle={An Example Article},
  pdfauthor={D. Doe, P. T. Frank, and J. E. Smith}
}
\fi

\usepackage{algorithmic}

\newsiamremark{assumption}{Assumption} % <- Preamble
\Crefname{ALC@unique}{Line}{Lines} % <- Preamble

\usepackage{booktabs}
\usepackage{multirow}

\begin{document}
\allowdisplaybreaks[4]
\maketitle

% REQUIRED
\begin{abstract}
Root mean square propagation (abbreviated as RMSProp) is a first-order stochastic algorithm used in machine learning widely. In this paper, a stable gradient-adjusted RMSProp (abbreviated as SGA-RMSProp) with mini-batch stochastic gradient is proposed, and its properties are studied on the linear least squares problem. R-linear convergence of the algorithm is established on the consistent linear least squares problem. The algorithm is also proved to converge R-linearly to a neighborhood of the minimizer for the inconsistent case, with the region of the neighborhood being controlled by the batch size. Furthermore, numerical experiments are conducted to compare the performances of SGA-RMSProp, stochastic gradient descent (abbreviated as SGD), and the original RMSProp with different batch sizes. The faster initial convergence rate of SGA-RMSProp is observed through numerical experiments and an adaptive strategy for switching from SGA-RMSProp to SGD is proposed, which combines the benefits of these two algorithms.
\end{abstract}

% REQUIRED

% REQUIRED

\section{Introduction}\label{sec: introduction}
To address the following finite-sum optimization problem
\begin{equation}
	\min_{\boldsymbol{x}\in \mathbb{R}^d}~ f(\boldsymbol{x}):=\sum_{i=1}^n f_i(\boldsymbol{x})
	\label{eq:optimization in machine learning}
\end{equation}
in machine learning, Tieleman and Hinton \cite{hinton2012lecture} developed a first-order stochastic optimization algorithm, called root mean square propagation (RMSProp). Its coordinate-wise version with time-varying hyperparameters is defined as follows:
\begin{align}
	&u_{k,j}=\beta_k u_{k-1, j}+\left(1-\beta_k\right) g_{k, j}^2, \label{eq: RMSProp u iter} \\
	&x_{k+1, j}=x_{k,j}-\eta_k g_{k,j} / \sqrt{u_{k,j}} \label{eq: RMSProp x iter}
\end{align}
for $j=1,\ldots,d$, where $x_{k,j}$ is the $j$th component of the $k$th iteration point $\boldsymbol{x}_k$ and $g_{k,j}$ denotes the $j$th component of the stochastic gradient $\boldsymbol{g}_k$. Due to the large data size $n$, the stochastic gradient is derived by sampling from $\{\nabla f_i(\boldsymbol{x}_k)\}_{i=1}^n$, which will be discussed later. The $j$th component of the moving average $\boldsymbol{u}_k$, denoted as $u_{k,j}$, is used to adjust the current stochastic gradient, which can be expressed as a linear combination of $u_{0,j}, g_{1,j}^2,\ldots, g_{k,j}^2$ as in \eqref{eq: RMSProp u iter}. The hyperparameters $\beta_k$ and $\eta_k$ are referred to as the discounting factor and the step size, respectively.

In this paper, we focus on the convergence rate of RMSProp \eqref{eq: RMSProp u iter}-\eqref{eq: RMSProp x iter} with the mini-batch stochastic gradient on the linear least squares problem (LLSP). It is well-known that LLSP has become one of the cornerstone problems in the optimization, machine learning, and statistics communities. In practice, it is widely used in computer vision \cite{naseem2010linear}, video signal processing \cite{zhang2019multiple}, calibration \cite{ling2018self}, linear classification \cite{cassioli2013incremental}, and financial markets \cite{lai2008statistical}.
Specifically, let $\boldsymbol{A}$ be an $n\times d$ full column rank real matrix. Then, LLSP can be formulated as
\begin{equation}\label{eq:LLSP obj}
	\min_{\boldsymbol{x}\in \mathbb{R}^d} \quad \frac{1}{2}\|\boldsymbol{A} \boldsymbol{x}-\boldsymbol{b}\|_2^2,
\end{equation}
where $\boldsymbol{x}\in \mathbb{R}^d$ represents the parameters to be estimated in the model, while  $\boldsymbol{A}=(\boldsymbol{a}_1,\ldots,\boldsymbol{a}_n)^{\top}$ and $\boldsymbol{b}=(b_1,\ldots,b_n)^{\top}$ are the observed data. Clearly, \eqref{eq:LLSP obj} can be viewed as a specific case of \eqref{eq:optimization in machine learning} with $f_i(\boldsymbol{x})=\frac{1}{2}(\boldsymbol{a}_i^{\top} \boldsymbol{x}-b_i)^2$, $i=1,\ldots, n$. Throughout this paper, the consistent (inconsistent) linear system represents $\boldsymbol{A}\boldsymbol{x}^*=\boldsymbol{b}$ ($\boldsymbol{A}\boldsymbol{x}^*\neq \boldsymbol{b}$) and the corresponding \eqref{eq:LLSP obj} is referred to as the consistent (inconsistent) LLSP, where $\boldsymbol{x}^*$ is the minimizer of \eqref{eq:LLSP obj}. 

At the $k$th iteration, let $\xi_{1}^{(k)},\ldots,$ $ \xi_{B}^{(k)}$ be a set of $B$ independent and identically distributed random variables, in which $\xi_{i}^{(k)}$ takes values in $\{1,\ldots,n\}$ with
$\mathbb{P}(\xi_{i}^{(k)}=j)=p_j>0$, $i=1,\ldots,B$, $j=1,\ldots,n$. Then the mini-batch stochastic gradient is given by
\begin{equation}\label{eq: k-th stochastic gradient}
	\boldsymbol{g}_k:=\frac{1}{B} \sum_{i=1}^B \frac{1}{p_{\scriptscriptstyle \xi_{i}^{(k)}}} \nabla f_{\scriptscriptstyle \xi_{i}^{(k)}}(\boldsymbol{x}_k)=\frac{1}{B} \sum_{i=1}^B \frac{1}{p_{\scriptscriptstyle \xi_{i}^{(k)}}} \boldsymbol{a}_{\scriptscriptstyle \xi_{i}^{(k)}}\left(\boldsymbol{a}_{\scriptscriptstyle \xi_{i}^{(k)}}^{\top} \boldsymbol{x}_k-b_{\scriptscriptstyle \xi_{i}^{(k)}}\right).
\end{equation}
In addition, aggregating \eqref{eq: RMSProp x iter} across all dimensions yields the vector form
\begin{equation}
	\boldsymbol{x}_{k+1}=\boldsymbol{x}_{k}-\eta_k \boldsymbol{D}_k\boldsymbol{g}_k,
	\label{eq: vector form of x in RMSProp}
\end{equation}
where $\boldsymbol{D}_k:=\text{diag}(\tfrac{1}{\sqrt{u_{k,1}}},\ldots,\tfrac{1}{\sqrt{u_{k,d}}})$. 

Equation \eqref{eq: vector form of x in RMSProp} reduces to the iteration scheme of stochastic gradient descent (SGD) when $\boldsymbol{D}_k\equiv I$. Strohmer and Vershynin \cite{strohmer2009randomized} proposed the randomized Kaczmarz algorithm and proved that it has a Q-linear convergence rate on consistent linear systems, while Needell et al. \cite{needell2016stochastic} showed that studying the randomized Kaczmarz algorithm on linear systems can be viewed as analyzing a special case of SGD on the associated LLSP. Notice that the objective function of LLSP in \eqref{eq:LLSP obj} is strongly convex and its gradient is Lipschitz continuous. Moulines and Bach \cite{moulines2011non} proved that SGD converges to a neighborhood of the minimizer R-linearly for such objective functions. Needell et al. \cite{needell2016stochastic} reduced the constant in this convergence rate.

In the above studies on SGD, $\eta_k$ is set to either a constant or a multiple of $1/k^\ell$ for $\ell\in[0,1]$. The Adagrad-norm algorithm studied by Ward et al. \cite{ward2020adagrad} uses a different step size and is more similar to RMSProp. Its iterative schemes are
\begin{align}
	&v_{k}=v_{k-1}+\|\boldsymbol{g}_{k}\|_2^2, \label{eq: Adagradnorm iter 1} \\
	&\boldsymbol{x}_{k+1}=\boldsymbol{x}_{k}-\frac{\eta_k}{v_{k}} \boldsymbol{g}_{k}, \label{eq: Adagradnorm iter 2}
\end{align}
where $v_k$ is a scalar. The distinction between \eqref{eq: Adagradnorm iter 1}-\eqref{eq: Adagradnorm iter 2} and RMSProp \eqref{eq: RMSProp u iter}-\eqref{eq: RMSProp x iter} is that the former still uses the stochastic gradient as the descent direction, that is, $\boldsymbol{D}_k\equiv\boldsymbol{I}$, while RMSProp modifies the descent direction. Xie et al. \cite{xie2020linear} proved that the Adagrad-norm algorithm achieves an R-linear convergence rate when the objective function is strongly convex or satisfies the Polyak-{\L}ojasiewicz inequality and $\nabla f_i\left(\boldsymbol{x}^*\right)=\boldsymbol{0}$, $i=1,\ldots,n$. Similar research includes \cite{wu2022adaloss}.

Moreover, \hspace{-0.2mm}some first-order\hspace{-0.2mm} preconditioned\hspace{-0.2mm} stochastic \hspace{-0.2mm}methods \hspace{-0.2mm}also have\hspace{-0.2mm} the iterative\hspace{-0.2mm} scheme of \eqref{eq: vector form of x in RMSProp}. For instance, Liu et al. \cite{liu2019acceleration} discussed how a fixed preconditioner $\boldsymbol{D}_k\equiv\boldsymbol{D}$ affects the performance of stochastic variance reduction gradient descent. The sequence $\{f(\boldsymbol{x}_k)\}$ generated by their algorithm was proved to converge R-linearly when each $f_i$ is strongly convex and has a Lipschitz continuous gradient. For the adaptive sketching-based preconditioners proposed by Lacotte and Pilanci \cite{lacotte2021fast}, the $\boldsymbol{g}_k$ and $\boldsymbol{D}_k$ in \eqref{eq: vector form of x in RMSProp} are chosen as $\nabla f(\boldsymbol{x}_k)$ and the inverse of $\boldsymbol{A}^{\top}\boldsymbol{L}^{\top}\boldsymbol{L}\boldsymbol{A}$ for LLSP \eqref{eq:LLSP obj} respectively, where $\boldsymbol{L}$ is a random matrix with entries independently drawn from specific distributions. They proved that the algorithm converges R-linearly. Besides, linear convergence rates are also attained by some other first-order stochastic algorithms, such as proximal SGD \cite{richtarik2020stochastic,gorbunov2020unified,khaled2023unified}, SGD with momentum \cite{loizou2020momentum,bollapragada2023fast,wang2023generalized,zeng2024adaptive}, randomized $r$-sets-Douglas-Rachford method \cite{han2024randomized}, and randomized primal-dual gradient method \cite{lan2018optimal}, though their iterative schemes are different from \eqref{eq: vector form of x in RMSProp}.

RMSProp constructs $\boldsymbol{D}_k$ from $\boldsymbol{g}_k$ and $\boldsymbol{D}_{k-1}$, making $\boldsymbol{D}_k$ a random and dynamically changing matrix that uses only first-order information. Liu et al. \cite{liu2022hyper} proved that, assuming a non-convex objective function and a bounded second-order moment of the infinity norm of the stochastic gradient, $\{\boldsymbol{x}_k\}$ generated by RMSProp satisfies $\frac{1}{K} \sum_{k=1}^K(\mathbb{E}[\|\nabla f(\boldsymbol{x}_k)\|_2^{4/3}])^{3/2}= \mathcal{O}(\frac{\log K}{\sqrt{K}})$ with $\eta_k\equiv\frac{1}{\sqrt{K}}$ and $\beta_k=\left(1-\frac{1}{k}\right)^{\ell}$, where $\ell$ is a positive constant and $K$ is the number of iterations. Under the assumptions that $f$ is gradient Lipschitz continuous and the stochastic gradient has coordinate-wise bounded noise variance, Li and Lin \cite{li2024frac} proved that $\frac{1}{K} \sum_{k=1}^K \mathbb{E}\left[\left\|\nabla f\left(\boldsymbol{x}_k\right)\right\|_1\right] \leq \mathcal{O}\left(\frac{\sqrt{d}}{K^{1/4}}\right)$ with $\eta_k\equiv\frac{\ell}{\sqrt{dK}}$ and $\beta_k\equiv1-\frac{1}{K}$, in which $\ell$ is a positive constant. Other studies have investigated the more general Adam-style algorithm, such as \cite{kingma2014adam,zou2019sufficient,guo2021novel,chen2022towards}, whose coordinate-wise version is defined as below:
\begin{align*}
	&m_{k,j}=\alpha_k m_{k-1, j}+\left(1-\alpha_k\right) g_{k, j},  \\
	&u_{k,j}=\beta_k u_{k-1, j}+\left(1-\beta_k\right) g_{k, j}^2, \\
	&x_{k+1, j}=x_{k,j}-\eta_k m_{k,j} / \sqrt{u_{k,j}} 
\end{align*}
for $j=1,\ldots,d$. Compared with RMSProp, the Adam-style algorithm uses an additional momentum $m_{k,j}$ to improve the performance of the algorithm. Furthermore, setting $\alpha_k=0$ causes the Adam-style algorithm to degenerate into RMSProp. Therefore, in the research on the convergence rate of Adam-style algorithm, allowing $\alpha_k=0$ makes the results applicable to  RMSProp. For instance, Zou et al. \cite{zou2019sufficient} and Chen et al. \cite{chen2022towards} studied the convergence rate of the Adam-style algorithm under different parameter settings. In their studies, the range of $\alpha_k$ always includes $0$, so the results can apply to  RMSProp. An aspect of these studies is that they dealt with general functions, and thus the resulting convergence rates are sublinear.

In numerical experiments, we observed that RMSProp, when implemented on LLSP, exhibits a linear convergence rate empirically under specific parameter settings. We further discovered that in numerous cases, RMSProp outperforms SGD in the early stages of iterations. Related observation of acceleration can also be noted in \cite{ma2022qualitative}. Since tuning the discounting factor $\beta_k$ can be cumbersome and sensitive, we introduce a stable gradient-adjusted mechanism that selects $\beta_k$ at each iteration. In this paper, we report the aforementioned experimental findings and explore whether RMSProp can attain a faster convergence rate on LLSP theoretically, with a focus on the linear convergence rate. It is worth mentioning that the convergence rate we obtain is better than sublinear, as we consider a specific LLSP.

The main contributions of this paper are summarized as follows:
\begin{itemize}
	\item We introduce a parameter $\varepsilon>0$, which is called adjusted level hereinafter, to stably control the adjustment applied to the stochastic gradient. Combining with a specific selection strategy of $\beta_k$, we guarantee $\|\boldsymbol{D}_k^{-1}\boldsymbol{D}_{k-1}-\boldsymbol{I}\|_2\leq \varepsilon$ during the iterations of RMSProp. So the algorithm is named stable gradient-adjusted RMSProp (SGA-RMSProp). Other properties arising from this stable gradient-adjusted mechanism are studied as well.
	
	\item We prove that SGA-RMSProp with the mini-batch stochastic gradient achieves an R-linear convergence rate on the consistent LLSP. For the inconsistent LLSP, we prove that SGA-RMSProp converges R-linearly to a neighborhood of the minimizer $\boldsymbol{x}^*$, where the region of the neighborhood is controlled by the batch size $B$. Moreover, we analytically discuss how the selection of the adjusted level $\varepsilon$ affects the performance of SGA-RMSProp.
	\item  We conduct numerical experiments to verify the R-linear convergence rate. In addition, the performances of SGA-RMSProp, SGD, and the original RMSProp (as defined by \eqref{eq: RMSProp u iter}-\eqref{eq: RMSProp x iter} with constant hyperparameters) on LLSP are compared numerically.
	\item We propose an adaptive condition for switching from SGA-RMSProp to SGD. Numerical experiments show that this method generally outperforms the vanilla SGD on LLSP in terms of computational time.
\end{itemize}

The remainder of this paper is organized as below. In \cref{sec:preliminaries}, notations and some results from linear algebra and probability theory are provided. In \cref{sec:SGA-RMSProp}, we formally introduce the framework of SGA-RMSProp and discuss its properties. In \cref{sec:convergence analysis of SGA-RMSProp}, we present our main results regarding the R-linear convergence rate of SGA-RMSProp. We report our numerical experiment results in \cref{sec:numerical experiments}. Finally, concluding remarks are presented in \cref{sec:concluding remarks}.

\section{Preliminaries}\label{sec:preliminaries}
This section introduces the notations and some established results serving as foundational components in our proofs.
\subsection{Notations}
Let $\mathbb{R}^{n}$, $\mathbb{R}^{n\times d}$, and $\mathbb{N}$ denote the sets of real $n$-dimensional vectors, $n \times d$ matrices, and natural numbers, respectively.  The sets of real $n\times n$ symmetric and positive semidefinite matrices are denoted by $\mathbb{S}^{n}$ and $\mathbb{S}^{n}_{+}$, respectively. For a vector $\boldsymbol{x}\in \mathbb{R}^d$, let $\|\boldsymbol{x}\|_2:=(\sum_{i=1}^{d}x_i^2)^{\frac{1}{2}}$, $\frac{1}{\sqrt{\boldsymbol{x}}}:=(\frac{1}{\sqrt{x_1}},\ldots,\frac{1}{\sqrt{x_d}})^{\top}$, and $\boldsymbol{x}^2:=\left(x_1^2,\ldots,x_d^2\right)^{\top}$. The Hadamard product is denoted by $\circ$. Given two matrices $\boldsymbol{X},\boldsymbol{Y} \in\mathbb{S}^{n}$, $\boldsymbol{X}\preceq\boldsymbol{Y}$ represents $\boldsymbol{Y}-\boldsymbol{X}\in\mathbb{S}^{n}_{+}$. The $j$th largest eigenvalue and the $j$th largest singular value of $\boldsymbol{X}$ are denoted by $\lambda_{j}$ and $s_{j}$ for $j=1,\ldots,d$, respectively. Let $\|\boldsymbol{X}\|_{\mathrm{F}}$ and $\|\boldsymbol{X}\|_2$ denote the Frobenius norm and the spectral norm which is equal to the maximum singular value of $\boldsymbol{X}$, respectively.  Let $\boldsymbol{X}=\text{diag}(x_{11},\ldots,x_{dd})$ be a diagonal matrix whose diagonal entries are $x_{11},\ldots,x_{dd}$. The matrix $\boldsymbol{X}^{\ell}$ denotes the diagonal matrix $\text{diag}(x_{11}^\ell,\ldots,x_{dd}^\ell)$ for a real number $\ell$. The symbols $\boldsymbol{0}$, $\boldsymbol{O}$, and $\boldsymbol{I}$ denote the zero column vector, zero matrix, and the identity matrix whose sizes vary from the context, respectively. The triple $(\Omega,\mathcal{A},\mathbb{P})$ denotes a probability space, where $\Omega$, $\mathcal{A}$, and $\mathbb{P}$ denote the sample space, $\sigma$-field of subsets of $\Omega$, and probability measure, respectively. For a random variable $\xi$ and a $\sigma$-field $\mathcal{F}\subset \mathcal{A}$, $\mathbb{E}[\cdot]$ and $\mathbb{E}[\xi|\mathcal{F}]$ denote the expectation and the conditional expectation of $\xi$ given $\mathcal{F}$, respectively.
In the remainder of this paper, $f(\boldsymbol{x})$ represents the object function $\frac{1}{2}\|\boldsymbol{A} \boldsymbol{x}-\boldsymbol{b}\|_2^2$ in \eqref{eq:LLSP obj}, while $f_i(\boldsymbol{x})$ denotes $\frac{1}{2}(\boldsymbol{a}_i^{\top} \boldsymbol{x}-b_i)^2, i=1, \ldots, n$.
\subsection{Linear algebra}
Next, two lemmas in linear algebra are presented.
\begin{lemma}[Weyl \cite{horn2012matrix}{~Theorem 4.3.1}]\label{lem:Weyl}
	Suppose that $\boldsymbol{X},\boldsymbol{Y}\in \mathbb{S}^{d}$. It holds that
	\begin{equation}\label{eq: Weyl inequality}
		\max _{i+k=j+d}\left\{\lambda_i(\boldsymbol{X})+\lambda_k(\boldsymbol{Y})\right\} \leq \lambda_j(\boldsymbol{X}+\boldsymbol{Y}) \leq \min _{i+k=j+1}\left\{\lambda_i(\boldsymbol{X})+\lambda_k(\boldsymbol{Y})\right\}
	\end{equation}
	for $1 \leq j \leq d$, where $\lambda_j$ represents the $j$th largest eigenvalue of the matrix.
\end{lemma}

The above result describes how the eigenvalues of a symmetric matrix $\boldsymbol{X}$ may change when a symmetric matrix $\boldsymbol{Y}$ is added \cite{horn2012matrix}. The next lemma is similar but deals with singular values.
\begin{lemma}[Ky Fan \cite{horn1994topics}{~Theorem 3.3.16}]	\label{lem:Ky Fan}
	Suppose that $\boldsymbol{X}, \boldsymbol{Y} \in \mathbb{R}^{d\times d}$. It holds that
	\begin{equation}\label{lem: Fan inequality}
		s_{i+j-1}(\boldsymbol{X}+\boldsymbol{Y})\leq s_i(\boldsymbol{X})+s_j(\boldsymbol{Y}), \quad s_{i+j-1}(\boldsymbol{X}\boldsymbol{Y})\leq s_i(\boldsymbol{X}) s_j(\boldsymbol{Y})
	\end{equation}
	for $ 1 \leq i, j \leq d, i+j-1 \leq d$, where $s_j$ represents the $j$th largest singular value of the matrix.
\end{lemma}

Here are some useful facts about eigenvalues and singular values of a matrix \cite{horn2012matrix}. Let $\boldsymbol{X}\in \mathbb{R}^{d\times d}$ be a real matrix, $s_j\left(\boldsymbol{X}\right)=\sqrt{\lambda_j\left(\boldsymbol{X}^{\top}\boldsymbol{X}\right)}$, $j=1,\ldots,d$, $\left\|\boldsymbol{X}\right\|_2=s_1\left(\boldsymbol{X}\right)=\sqrt{\lambda_1\left(\boldsymbol{X}^{\top}\boldsymbol{X}\right)}$. If $\boldsymbol{X}$ is positive semidefinite, then $\lambda_j(\boldsymbol{X})=s_j(\boldsymbol{X})$, $j=1,\ldots,d$. Moreover, if $\boldsymbol{X}=\text{diag}(x_{11},\ldots,x_{dd})$ is a diagonal matrix, then $\left\|\boldsymbol{X}\right\|_2=s_1\left(\boldsymbol{X}\right)=\max\limits_{j=1,\ldots,d} \left|x_{jj}\right|$.
\subsection{Probability theory}
Taking $(\Omega,\mathcal{A},\mathbb{P})$ as the probability space, we review some notions and results from probability theory. 

\begin{definition}[Filtration \cite{durrett2019probability}]\label{def:filtration} Let $\mathcal{K}\subset \mathbb{N}$ be an index set. A filtration $\{\mathcal{F}_k\}_{k\in\mathcal{K}}$ is an increasing sequence of $\sigma$-fields, that is, $\mathcal{F}_{k_1}\subset\mathcal{F}_{k_2}\subset\mathcal{A}$ for all $k_1, k_2\in\mathcal{K}$ with $k_1\leq k_2$.
\end{definition}

Consider $\xi$, $\zeta$ as $\mathcal{A}$-measurable random variables with $\mathbb{E}[|\xi|]$, $\mathbb{E}[|\zeta|]$, $\mathbb{E}[|\xi\zeta|] < +\infty$ and $\mathcal{F}\subset \mathcal{A}$ as a $\sigma$-field. Since the $\sigma$-fields we consider in this paper have no null sets other than the empty set, we omit the ``almost surely'' in equalities like $\xi=\mathbb{E}[\zeta|\mathcal{F}]$ throughout this paper. The following are three useful propositions.
\begin{proposition}[\cite{durrett2019probability}]
	If $\xi$ is $\mathcal{F}$-measurable, then
	\begin{equation}\label{eq: CE measurable}
		\mathbb{E}[\xi \zeta|\mathcal{F}]=\xi \mathbb{E}[\zeta|\mathcal{F}],~\text{and specifically,}~\mathbb{E}[\xi|\mathcal{F}]=\xi.
	\end{equation}
\end{proposition}
\begin{proposition}[\cite{durrett2019probability}]
	If the $\sigma$-field generated by $\xi$ is independent of $\mathcal{F}$, which is also referred to as $\xi$ is independent of $\mathcal{F}$, then
	\begin{equation}\label{eq: CE independent}
		\mathbb{E}[\xi|\mathcal{F}]=\mathbb{E}[\xi].
	\end{equation}
\end{proposition}
\begin{proposition}[\cite{durrett2019probability}]\label{prop:CE tower}
	If $\mathcal{G}$ is a $\sigma$-field with $\mathcal{G} \subset \mathcal{F}$, then
\begin{equation}\label{eq: CE tower}
		\mathbb{E}[\mathbb{E}[\xi|\mathcal{F}]|\mathcal{G}]=\mathbb{E}[\xi|\mathcal{G}].
\end{equation}
\end{proposition}

To estimate the expectation of the spectral norm of independent, zero mean random matrices sum, we introduce the matrix Bernstein inequality.
\begin{lemma}[Matrix Bernstein \cite{tropp2015an}{~Theorem 6.1.1}] \label{lem:lemma Matrix Bernstein}
	Let $\left\{\boldsymbol{Z}_i\right\}_{i=1}^{n}$ be independent, random $d_1 \times d_2$ matrices, $\mathbb{E}[\boldsymbol{Z}_i]=\boldsymbol{O}$, $\left\|\boldsymbol{Z}_i\right\|_2 \leq L$, $\boldsymbol{Z}:=\sum_{i=1}^{n} \boldsymbol{Z}_i$, and
	\begin{equation}\label{eq: matrix variance statistic}
		\nu(\boldsymbol{Z}) =\max \left\{\left\|\sum_{i=1}^n \mathbb{E}[\boldsymbol{Z}_i \boldsymbol{Z}_i^{\top}]\right\|_2,\left\|\sum_{i=1}^n \mathbb{E}[\boldsymbol{Z}_i^{\top} \boldsymbol{Z}_i]\right\|_2\right\}
	\end{equation}
be the matrix variance statistic of the sum.	
Then we have
$$\mathbb{E}[\|\boldsymbol{Z}\|_2] \leq \sqrt{2 v(\boldsymbol{Z}) \log(d_1+d_2)}+\frac{L \log \left(d_1+d_2\right)}{3}.$$
\end{lemma}

\section{SGA-RMSProp}\label{sec:SGA-RMSProp}
In this section, we first present the pseudocode of \hyperref[alg:SGA-RMSProp]{SGA-RMSProp} and then discuss some properties of the algorithm.
\begin{algorithm}[H]
	\caption{SGA-RMSProp}\label{alg:SGA-RMSProp}
	\begin{algorithmic}[1]
		\REQUIRE Step size $\{\eta_k\}$, adjusted level $\varepsilon$, lower and upper bounds $\underline{u}$, $\overline{u}$, initial value $\boldsymbol{x}_1$.
		\STATE $\boldsymbol{u}_{0}=\left(\frac{1}{\overline{u}^2},\ldots,\frac{1}{\overline{u}^2}\right)^{\top}$.
		\FOR{$k=1,\ldots,K$}
		\STATE Sample a stochastic gradient $\boldsymbol{g}_k$.
		\STATE $\beta_k=\text{ \textbf{$\boldsymbol{\beta}$-Selection}}(\boldsymbol{g}_k,\boldsymbol{u_{k-1}},\varepsilon,\underline{u},\overline{u})$.\label{line: selection in SGA-RMSProp}
		\STATE $\boldsymbol{u}_k=\beta_k \boldsymbol{u}_{k-1}+\left(1-\beta_k\right) \boldsymbol{g}_k^2$.
		\STATE $\boldsymbol{x}_{k+1}=\boldsymbol{x}_k-\frac{\eta_k}{\sqrt{\boldsymbol{u}_k}} \circ \boldsymbol{g}_k$.
		\ENDFOR
		\ENSURE The optimal point $\boldsymbol{x}_{K+1}$.
	\end{algorithmic}
\end{algorithm}
The algorithm above can be viewed as a vector form of RMSProp \eqref{eq: RMSProp u iter}-\eqref{eq: RMSProp x iter}, in which the discounting factor $\beta_k$ is selected during each iteration, as shown in \Cref{line: selection in SGA-RMSProp}. The details of \hyperref[alg:beta selection]{\textbf{$\beta$-Selection}} are presented below, and its validity will be illustrated through \cref{lem:SGA-RMSProp control}.

\begin{algorithm}[H]\label{alg:beta selection}
	\caption{$\beta$-Selection$(\boldsymbol{g}_k,\boldsymbol{u}_{k-1},\varepsilon,\underline{u},\overline{u})$}
	\begin{algorithmic}[1]
		\REQUIRE Stochastic gradient $\boldsymbol{g}_{k}$, moving average $\boldsymbol{u}_{k-1}$,  adjusted level $\varepsilon$, lower and upper bounds $\underline{u}$, $\overline{u}$.
		\IF{$\exists~j \in\{1,\ldots,d\}$ such that $\frac{g_{k,j}^2}{u_{k-1,j}}\leq1$}
		\STATE $\beta_k= 1$.
		\ELSE
		\STATE Select $\beta_k\in\left[\max\limits_{j=1, \ldots, d}\left\{\frac{g_{k, j}^2-\frac{1}{\underline{u}^2}}{g_{k, j}^2-u_{k-1, j}}, \frac{g_{k, j}^2-(1+\varepsilon)^2 u_{k-1, j}}{g_{k, j}^2-u_{k-1, j}},0\right\},1\right]$.
		\ENDIF
		\ENSURE The discounting factor $\beta_k$.
	\end{algorithmic}
\end{algorithm}

Recall from \eqref{eq: k-th stochastic gradient} that the mini-batch stochastic gradient in the $k$th iteration is given by
\begin{equation*}
	\boldsymbol{g}_k:=\frac{1}{B} \sum_{i=1}^B \frac{1}{p_{\scriptscriptstyle \xi_{i}^{(k)}}} \nabla f_{\scriptscriptstyle \xi_{i}^{(k)}}(\boldsymbol{x}_k)=\frac{1}{B} \sum_{i=1}^B \frac{1}{p_{\scriptscriptstyle \xi_{i}^{(k)}}} \boldsymbol{a}_{\scriptscriptstyle \xi_{i}^{(k)}}\left(\boldsymbol{a}_{\scriptscriptstyle \xi_{i}^{(k)}}^{\top} \boldsymbol{x}_k-b_{\scriptscriptstyle \xi_{i}^{(k)}}\right).
\end{equation*}
Let $\mathcal{S}_k$ denote the set $\{\xi_{1}^{(k)},\ldots,$ $ \xi_{B}^{(k)}\}$. It can be easily verified that the randomness of $\boldsymbol{g}_k$ arises not only from $\mathcal{S}_k$ but also from $\boldsymbol{x}_k$ for each $k=2,\ldots,K$.
Let $\mathcal{F}_k$ be the $\sigma$-field generated by the random variables in $\mathcal{S}_1,\ldots,\mathcal{S}_{k-1}$, and then $\{\mathcal{F}_k\}_{k=2}^K$ is a filtration by \cref{def:filtration}. Since $\boldsymbol{x}_k$ is $\mathcal{F}_k$-measurable and $\xi_{i}^{(k)}\in \mathcal{S}_k$ is independent of $\mathcal{F}_k$, $i=1,\ldots,B$, we have
{\small
\begin{align}\label{eq: unbiased stochastic gradient}
		\mathbb{E}\left[\frac{1}{p_{\scriptscriptstyle \xi_{i}^{(k)}}} \boldsymbol{a}_{\scriptscriptstyle \xi_{i}^{(k)}}\left(\boldsymbol{a}_{\scriptscriptstyle \xi_{i}^{(k)}}^{\top} \boldsymbol{x}_k-b_{\scriptscriptstyle \xi_{i}^{(k)}}\right)\Bigg|\mathcal{F}_k\right]&\!\overset{\eqref{eq: CE measurable}}{=}\! \mathbb{E}\!\left[\frac{1}{p_{\scriptscriptstyle \xi_{i}^{(k)}}} \boldsymbol{a}_{\scriptscriptstyle \xi_{i}^{(k)}}\boldsymbol{a}_{\scriptscriptstyle \xi_{i}^{(k)}}^{\top} \Bigg|\mathcal{F}_k\right]\boldsymbol{x}_k-\mathbb{E}\!\left[\frac{1}{p_{\scriptscriptstyle \xi_{i}^{(k)}}} \boldsymbol{a}_{\scriptscriptstyle \xi_{i}^{(k)}}b_{\scriptscriptstyle \xi_{i}^{(k)}}\Bigg|\mathcal{F}_k\right] \nonumber\\
	&\!\overset{\eqref{eq: CE independent}}{=}\!\mathbb{E}\!\left[\frac{1}{p_{\scriptscriptstyle \xi_{i}^{(k)}}} \boldsymbol{a}_{\scriptscriptstyle \xi_{i}^{(k)}}\boldsymbol{a}_{\scriptscriptstyle \xi_{i}^{(k)}}^{\top}\right]\boldsymbol{x}_k\!-\!\mathbb{E}\!\left[\frac{1}{p_{\scriptscriptstyle \xi_{i}^{(k)}}} \boldsymbol{a}_{\scriptscriptstyle \xi_{i}^{(k)}}b_{\scriptscriptstyle \xi_{i}^{(k)}}\right] \!=\!\sum_{i=1}^{n}\boldsymbol{a}_{i}\left(\boldsymbol{a}_{i}^{\top} \boldsymbol{x}_k-b_{i}\right)
\end{align}
}
\hspace{-2mm}for each $i=1,\ldots,B$, and then $\mathbb{E}[\boldsymbol{g}_k|\mathcal{F}_k]=\nabla f(\boldsymbol{x}_k)$, $k=2,\ldots,K$. Similarly, we have $\mathbb{E}[\boldsymbol{g}_1]=\nabla f(\boldsymbol{x}_1)$. So the stochastic gradient is unbiased. 

Next, recalling that $\boldsymbol{D}_k=\text{diag}(\frac{1}{\sqrt{u_{k,1}}},\ldots,\frac{1}{\sqrt{u_{k,d}}})$, $k=0,\ldots,K$, we present the following proposition to illustrate how \hyperref[alg:SGA-RMSProp]{SGA-RMSProp} uses $\varepsilon$ to modulate the variation of the adjustment made to the stochastic gradient, thereby achieving stable control.
\begin{proposition} \label{lem:SGA-RMSProp control}Let $\{\boldsymbol{u}_{k}\}_{k=0}^{K}$ be the sequence generated by \hyperref[alg:SGA-RMSProp]{SGA-RMSProp}. Then $u_{k,j}$ satisfies that $\underline{u} \leq \frac{1}{\sqrt{u_{k, j}}} \leq \overline{u}$, $k=0,\ldots,K$, $j=1, \ldots, d$. Equivalently, for each $k=0,\ldots,K$,
	\vspace{-4mm}
	\begin{equation}\label{eq: D_k estimate}
		\left\|\boldsymbol{D}_k^\ell\right\|_2\leq \overline{u}^\ell,\quad \left\|\boldsymbol{D}_k^{-\ell}\right\|_2 \leq \underline{u}^{-\ell},
	\end{equation}
	where $\ell$ is a positive number.
	Moreover, the matrix $\boldsymbol{D}_k$ satisfies
	\begin{equation}\label{eq:D_k D_k-1 estimate}
		\left\|\boldsymbol{D}_k^{\frac{1}{2}} \boldsymbol{D}_{k-1}^{-\frac{1}{2}}\right\|_2  \leq 1 ~\text{and}~ \left\|\boldsymbol{D}_k^{-1} \boldsymbol{D}_{k-1}-\boldsymbol{I}\right\|_2 \leq \varepsilon,  k=1,\ldots,K.
	\end{equation}
\end{proposition}
\begin{proof}
	To prove the first part of this proposition, we employ induction on the variable $k$.
	When $k=0$, recalling the initial values of $u_{0,j}$ in \hyperref[alg:SGA-RMSProp]{SGA-RMSProp}, we have $\underline{u} \leq \frac{1}{\sqrt{u_{0, j}}} = \overline{u}$, $j=1, \ldots, d$. Next assume that $\underline{u} \leq \frac{1}{\sqrt{u_{k, j}}} \leq \overline{u}$ for each $k=0, \ldots, t$, $j =1, \ldots, d$.
	Now we show that $\underline{u} \leq \frac{1}{\sqrt{u_{t+1, j}}} \leq \overline{u}$, $j =1, \ldots, d$. Under the induction hypothesis, $\beta_{t+1} \in [0,1]$. If there exists $j$ such that $\frac{g_{t+1,j}^2}{u_{t,j}}\leq1$, then we have $\boldsymbol{u}_{t+1}=\boldsymbol{u}_t$ since $\beta_{t+1}=1$ and $\boldsymbol{u}_{t+1}=\beta_{t+1} \boldsymbol{u}_{t}+\left(1-\beta_{t+1}\right) \boldsymbol{g}_{t+1}^2$. In this case, we have $\underline{u} \leq \frac{1}{\sqrt{u_{t+1, j}}} \leq \overline{u}$. If $\frac{g_{t+1,j}^2}{u_{t,j}}> 1$ for each $j=1,\ldots, d$, then we have
	\begin{equation*}
		\beta_{t+1}\geq \max _{j=1, \ldots, d}\left\{\frac{g_{t+1, j}^2-\frac{1}{\underline{u}^2}}{g_{t+1, j}^2-u_{t, j}}, \frac{g_{t+1, j}^2-(1+\varepsilon)^2 u_{t, j}}{g_{t+1, j}^2-u_{t, j}},0\right\}\geq \frac{g_{t+1, j}^2-\frac{1}{\underline{u}^2}}{g_{t+1, j}^2-u_{t, j}},~j=1,\ldots,d,
	\end{equation*}
	which implies that $u_{t+1,j}=\beta_{t+1} u_{t,j}+\left(1-\beta_{t+1}\right) g_{t+1,j}^2\leq \frac{1}{\underline{u}^2}$, $j=1,\ldots,d$. In addition, since
	\begin{equation*}
		\frac{u_{t+1, j}}{u_{t, j}}=\frac{\beta_{t+1} u_{t, j}+\left(1-\beta_{t+1}\right) g_{t+1, j}^2}{u_{t, j}}=\beta_{t+1}+\left(1-\beta_{t+1}\right) \frac{g_{t+1, j}^2}{u_{t, j}}\geq 1, j=1,\ldots,d,
		\label{utdividedut1}
	\end{equation*}
	we have $u_{t+1,j}\geq u_{t,j}\geq \frac{1}{\overline{u}^2}$, $j=1,\ldots,d$. Hence,  $\underline{u} \leq \frac{1}{\sqrt{u_{k, j}}} \leq \overline{u}$ for each $j=1, \ldots, d$, $k= t+1$. Since $\boldsymbol{D}_k$ is a diagonal matrix with all its diagonal elements positive, it is invertible and  $\|\boldsymbol{D}_k^\ell\|_2\leq \overline{u}^\ell$, $\|\boldsymbol{D}_k^{-\ell}\|_2 \leq \underline{u}^{-\ell}$, $k=0,\ldots,K$.
	
	Next, we prove the second part of the lemma.  The proof of the first part yields $u_{k,j}\geq u_{k-1,j}$, $j=1,\ldots, d$, $k=1,\ldots,K$. Combined with $\|\boldsymbol{D}_k^{\frac{1}{2}} \boldsymbol{D}_{k-1}^{-\frac{1}{2}}\|_2=\max\limits_{j=1,\ldots,d} \{(\frac{u_{k-1,j}}{u_{k,j}})^{\frac{1}{4}}\}$, we conclude that $\|\boldsymbol{D}_k^{\frac{1}{2}} \boldsymbol{D}_{k-1}^{-\frac{1}{2}}\|_2  \leq 1$, $k=1,\ldots,K$.
	
	For $\|\boldsymbol{D}_k^{-1} \boldsymbol{D}_{k-1}-\boldsymbol{I}\|_2$, $k=1,\ldots,K$, we have
	\begin{equation}\label{eq: pace control}
		\left\|\boldsymbol{D}_k^{-1} \boldsymbol{D}_{k-1}-\boldsymbol{I}\right\|_2 \leq \varepsilon \Longleftrightarrow\left|\sqrt{\frac{u_{k, j}}{u_{k-1, j}}}-1\right| \leq \varepsilon \Longleftrightarrow 1-\varepsilon \leq \sqrt{\frac{u_{k, j}}{u_{k-1, j}}} \leq 1+\varepsilon
	\end{equation}
	for $j=1, \ldots, d$. Since $\frac{u_{k, j}}{u_{k-1, j}}\geq 1$ for each $k=1,\ldots,K$, $j=1,\ldots,d$, it suffices to consider solely the second inequality. By the following transformation,
	\begin{equation*}
		\frac{u_{k, j}}{u_{k-1, j}} \leq(1+\varepsilon)^2  \Longleftrightarrow\left(1-\frac{g_{k, j}^2}{u_{k-1, j}}\right) \beta_k+\frac{g_{k, j}^2}{u_{k-1, j}} \leq(1+\varepsilon)^2, \Longleftrightarrow \beta_k \geq \frac{g_{k, j}^2-(1+\varepsilon)^2 u_{k-1, j}}{g_{k, j}^2-u_{k-1, j}}
	\end{equation*}
	for $j=1,\ldots, d,~k=1,\ldots,K$, and recalling the selection of $\beta_k$, we have
	\begin{equation*}
		\beta_k \geq \max _{j=1,\ldots, d}\left\{\frac{g_{k, j}^2-\frac{1}{\underline{u}^2}}{g_{k, j}^2-u_{k-1, j}}, \frac{g_{k, j}^2-(1+\varepsilon)^2 u_{k-1, j}}{g_{k, j}^2-u_{k-1, j}},0\right\}\geq \frac{g_{k, j}^2-(1+\varepsilon)^2 u_{k-1, j}}{g_{k, j}^2-u_{k-1, j}},~k=1,\ldots,K.
	\end{equation*}
	
	The proof is completed.
\end{proof}

This proposition shows two aspects of \hyperref[alg:SGA-RMSProp]{SGA-RMSProp}. Firstly, our algorithm can control the range of all $u_{k,j}$. As a result, reasonable settings of $\underline{u}$ and $\overline{u}$ can avoid potential numerical issues caused by the division of $\frac{1}{\sqrt{u_{k,j}}}$. Secondly, the adjusted level $\varepsilon$ stably controls the impact of squared gradients on the adjustment applied to the stochastic gradient by $\|\boldsymbol{D}_k^{-1} \boldsymbol{D}_{k-1}-\boldsymbol{I}\|_2 \leq \varepsilon$. It provides us the ability to quantitatively manage the change of $u_{k,j}$ at each iteration. A larger $\varepsilon$ allows $\boldsymbol{D}_k$ to adjust more, while a smaller $\varepsilon$ has the opposite effect. 
 
\section{Convergence analysis of SGA-RMSProp}\label{sec:convergence analysis of SGA-RMSProp}
In this section, we first estimate the convergence rate of \hyperref[alg:SGA-RMSProp]{SGA-RMSProp} on the consistent LLSP, then on the inconsistent LLSP.
\subsection{Consistent LLSP}
Denote the minimizer of \eqref{eq:LLSP obj} as $\boldsymbol{x}^*$ , and then the consistency of  LLSP implies $\boldsymbol{A}\boldsymbol{x}^*=\boldsymbol{b}$, that is, $\boldsymbol{a}_i^{\top} \boldsymbol{x}^*=b_i$, $i=1,\ldots,n$. Replacing $b_i$ with $\boldsymbol{a}_i^{\top} \boldsymbol{x}^*$ in \eqref{eq: k-th stochastic gradient}, by \eqref{eq: vector form of x in RMSProp}, we have
\vspace{-5mm}
\begin{align*}
	\boldsymbol{x}_{k+1}-\boldsymbol{x}^* & =\boldsymbol{x}_k-\boldsymbol{x}^*-\frac{\eta_k}{B}\boldsymbol{D}_k \left( \sum_{i=1}^{B} \frac{1}{p_{\scriptscriptstyle \xi_{i}^{(k)}}} \boldsymbol{a}_{\scriptscriptstyle \xi_{i}^{(k)}} \boldsymbol{a}_{\scriptscriptstyle \xi_{i}^{(k)}}^{\top}\left(\boldsymbol{x}_k-\boldsymbol{x}^*\right) \right)  \\
	& =\boldsymbol{x}_k-\boldsymbol{x}^*-\eta_k \boldsymbol{D}_k \boldsymbol{M}_k\left(\boldsymbol{x}_k-\boldsymbol{x}^*\right)  \\
	& =\left(\boldsymbol{I}-\eta_k \boldsymbol{D}_k \boldsymbol{M}_k\right)\left(\boldsymbol{x}_k-\boldsymbol{x}^*\right),
\end{align*}
where $\boldsymbol{M}_k:=\frac{1}{B} \sum_{i=1}^{B} \frac{1}{p_{\scriptscriptstyle \xi_{i}^{(k)}}} \boldsymbol{a}_{\scriptscriptstyle \xi_{i}^{(k)}} \boldsymbol{a}_{\scriptscriptstyle \xi_{i}^{(k)}}^{\top}$ and $\boldsymbol{D}_k=\text{diag}(\frac{1}{\sqrt{u_{k,1}}},\ldots,\frac{1}{\sqrt{u_{k,d}}})$.

Denote $\boldsymbol{Y}_k: =\boldsymbol{I}-\eta_k \boldsymbol{D}_k \boldsymbol{M}_k$ as the stochastic transition matrix in the $k$th iteration. We have $\boldsymbol{x}_{K+1}-\boldsymbol{x}^*=\boldsymbol{Y}_K\left(\boldsymbol{x}_K-\boldsymbol{x}^*\right)=\ldots=\boldsymbol{Y}_K \cdots \boldsymbol{Y}_1\left(\boldsymbol{x}_1-\boldsymbol{x}^*\right)$
, which implies that
\begin{equation}\label{eq: norm inequality}
	\left\|\boldsymbol{x}_{K+1}-\boldsymbol{x}^*\right\|_2 \leq\left\|\boldsymbol{Y}_K  \cdots \boldsymbol{Y}_1\right\|_2\left\|\boldsymbol{x}_1-\boldsymbol{x}^*\right\|_2.
\end{equation}
Therefore, establishing the convergence rate of \hyperref[alg:SGA-RMSProp]{SGA-RMSProp} hinges upon the estimation of $\mathbb{E}\left[\left\|\boldsymbol{Y}_K  \cdots \boldsymbol{Y}_1\right\|_2\right]$. To achieve this, we present the following two lemmas. The first one shows the relationship between $p_j$ and $\boldsymbol{a}_j$, $j=1,\ldots,n$.
\begin{lemma}\label{lem:probability bound}
	Let $\boldsymbol{A}=(\boldsymbol{a}_1,\ldots,\boldsymbol{a}_n)^{\top}\in\mathbb{R}^{n\times d}$ be a non-zero matrix and $\{p_j\}_{j=1}^{n}$ be the non-zero probabilities, that is, $p_j>0$, $j=1,\ldots,n$, and $\sum_{j=1}^{n} p_j=1$. Then we have 
	\begin{equation*}
		\max_{j=1,\ldots,n}\left\{\frac{\|\boldsymbol{a}_j\|_{2}^{2}}{p_j}\right\} \geq \|\boldsymbol{A}\|_{\mathrm{F}}^2.
	\end{equation*}
\end{lemma}
\begin{proof}
	Assume that $\frac{\|\boldsymbol{a}_j\|_{2}^{2}}{p_j}<\|\boldsymbol{A}\|_{\mathrm{F}}^2$, $j=1,\ldots,n$. It follows that $\sum\limits_{j=1}^{n} p_j > \sum\limits_{j=1}^{n} \frac{\|\boldsymbol{a}_j\|_{2}^{2}}{\|\boldsymbol{A}\|_{\mathrm{F}}^2} =1$, a contradiction. Therefore, we have $\max\limits_{j=1,\ldots,n}\left\{\|\boldsymbol{a}_j\|_{2}^{2}/p_j\right\} \geq \|\boldsymbol{A}\|_{\mathrm{F}}^2$.
\end{proof}

Next, we give an upper bound on $\mathbb{E}[\|\boldsymbol{M}_k-\mathbb{E}[\boldsymbol{M}_k]\|_2]$, which is a monotonically decreasing function of the batch size, similar to Lemma 2 in \cite{bollapragada2023fast}. 
\begin{lemma} \label{lem: standard error upper bound}
	Let $\boldsymbol{M}_k=\frac{1}{B} \sum_{i=1}^{B} \frac{1}{p_{\scriptscriptstyle \xi_{i}^{(k)}}} \boldsymbol{a}_{\scriptscriptstyle \xi_{i}^{(k)}} \boldsymbol{a}_{\scriptscriptstyle \xi_{i}^{(k)}}^{\top}$, where $\xi_{i}^{(k)}$, $i=1,\ldots,B$, $k=1,\ldots,K$, are independent and identically distributed random variables, taking values in $\{1,\ldots,n\}$ with
	$\mathbb{P}(\xi_{i}^{(k)}=j)=p_j$, $j=1,\ldots,n$. Then, for each $k=1,\ldots,K$, we have 
	\begin{align} \label{eq: standard error upper bound}\textstyle
		\mathbb{E}[\|\boldsymbol{M}_k-\mathbb{E}[\boldsymbol{M}_k]\|_2]  \leq&\Bigg(\frac{2\Big(\max\limits_{j=1,\ldots,n}\left\{\frac{\|\boldsymbol{a}_j\|_{2}^{2}}{p_j}\right\}-\lambda_d\left(\boldsymbol{A}^{\top} \boldsymbol{A}\right)\Big)\|\boldsymbol{A}\|_2^2 \log (2 d)}{B}\Bigg)^{\frac{1}{2}} \nonumber \\
		&+ \frac{\log(2d)}{3B}\left(\max\limits_{j=1,\ldots,n}\left\{\frac{\|\boldsymbol{a}_j\|_{2}^{2}}{p_j}\right\}-\lambda_d\left(\boldsymbol{A}^{\top} \boldsymbol{A}\right)\right).
	\end{align}
\end{lemma}

\begin{proof}
	Since $\xi_{i}^{(k)}$, $k=1,\ldots,K$, are independent and identically distributed random variables,  
	our proof below holds for each $k=1,\ldots,K$. Similar to \eqref{eq: unbiased stochastic gradient}, we have $\mathbb{E}[\boldsymbol{M}_k]=\boldsymbol{A}^{\top}\boldsymbol{A}$. Denoting $\boldsymbol{Z}_{\scriptscriptstyle \xi_i^{(k)}}:=\frac{1}{B}\left(\frac{1}{{\scriptstyle p}_{\scriptscriptstyle \xi_{i}^{(k)}}} \boldsymbol{a}_{\scriptscriptstyle \xi_i^{(k)}} \boldsymbol{a}_{\scriptscriptstyle \xi_i^{(k)}}^{\top}-\boldsymbol{A}^{\top} \boldsymbol{A}\right)$ and $\boldsymbol{Z}_k:=\sum_{i=1}^{B} \boldsymbol{Z}_{\scriptscriptstyle \xi_i^{(k)}}$, it is clear that $\mathbb{E}[\|\boldsymbol{M}_k-\mathbb{E}[\boldsymbol{M}_k]\|_2]=\mathbb{E}(\|\boldsymbol{Z}_k\|_2)$ and $\mathbb{E}[\boldsymbol{Z}_{\scriptscriptstyle \xi_i^{(k)}}]=\boldsymbol{O}$. Since $\boldsymbol{Z}_{\scriptscriptstyle \xi_i^{(k)}}$ is a real symmetric matrix, $\|\boldsymbol{Z}_{\scriptscriptstyle \xi_i^{(k)}}\|_2=s_1(\boldsymbol{Z}_{\scriptscriptstyle \xi_i^{(k)}})=\max \{\lambda_1(\boldsymbol{Z}_{\scriptscriptstyle \xi_i^{(k)}}),-\lambda_d(\boldsymbol{Z}_{\scriptscriptstyle \xi_i^{(k)}})\}$, $i=1,\ldots,B$. By \cref{lem:Weyl}, we have
	\begin{equation*}
			\lambda_d\!\left(\!\boldsymbol{Z}_{\scriptscriptstyle \xi_i^{(k)}}\!\!\right)\!=\!\lambda_d\!\left(\!\frac{1}{B}\!\left(\!\frac{1}{p_{\scriptscriptstyle \xi_i^{(k)}\!}} \boldsymbol{a}_{\scriptscriptstyle\xi_i^{(k)}} \boldsymbol{a}_{\scriptscriptstyle \xi_i^{(k)}}^{\top}\!-\!\boldsymbol{A}^{\top} \!\boldsymbol{A}\!\right)\hspace{-0.4em}\right) \!\!\overset{\eqref{eq: Weyl inequality}}{\geq}\!\!\frac{1}{B}\!\!\left(\!\lambda_d\!\left(\!\frac{1}{p_{\scriptscriptstyle \xi_i^{(k)}}\!} \boldsymbol{a}_{\scriptscriptstyle \xi_i^{(k)}} \boldsymbol{a}_{\scriptscriptstyle \xi_i^{(k)}}^{\top}\!\right)\!-\!\lambda_1\!\left(\!\boldsymbol{A}^{\top} \!\boldsymbol{A}\!\right) \!\!\right)\! \geq \!-\frac{\lambda_1 \!\left(\!\boldsymbol{A}^{\top}\! \boldsymbol{A}\!\right)\!}{B},
	\end{equation*}
	and
	\begin{equation*}
			\lambda_1\left(\boldsymbol{Z}_{\scriptscriptstyle \xi_i^{(k)}}\right)\overset{\eqref{eq: Weyl inequality}}{\leq} \frac{1}{B}\left( \lambda_1\left(\frac{1}{p_{\scriptscriptstyle \xi_i^{(k)}}\!} \boldsymbol{a}_{\scriptscriptstyle \xi_i^{(k)}} \boldsymbol{a}_{\scriptscriptstyle \xi_i^{(k)}}^{\top}\right)- \lambda_d\left(\boldsymbol{A}^{\top} \boldsymbol{A}\right)\right) \leq\frac{\max\limits_{j=1,\ldots,n}\left\{\frac{\|\boldsymbol{a}_j\|_{2}^{2}}{p_j}\right\}-\lambda_d\left(\boldsymbol{A}^{\top} \boldsymbol{A}\right)}{B} .
	\end{equation*}
	Therefore, $\|\boldsymbol{Z}_{\scriptscriptstyle \xi_i^{(k)}}\|_2 \leq \max \{\frac{\lambda_1\left(\boldsymbol{A}^{\top} \boldsymbol{A}\right)}{B} , \frac{\max\limits_{j=1,\ldots,n}\left\{\|\boldsymbol{a}_j\|_{2}^{2}/p_j\right\}-\lambda_d\left(\boldsymbol{A}^{\top} \boldsymbol{A}\right)}{B}\}$. By the definition of Frobenius norm, we have  $\|\boldsymbol{A}\|_{\mathrm{F}}^2=\sum_{i=1}^d s_i^2(\boldsymbol{A})=\sum_{i=1}^d \lambda_i(\boldsymbol{A}^{\top} \boldsymbol{A})$, which implies $\|\boldsymbol{A}\|_{\mathrm{F}}^2 \geq \lambda_1(\boldsymbol{A}^{\top} \boldsymbol{A})+\lambda_d(\boldsymbol{A}^{\top} \boldsymbol{A})$. By \cref{lem:probability bound}, it follows that $\max\limits_{j=1,\ldots,n}\{\|\boldsymbol{a}_j\|_{2}^{2}/p_j\}-\lambda_d(\boldsymbol{A}^{\top} \boldsymbol{A})\geq\lambda_1(\boldsymbol{A}^{\top} \boldsymbol{A})$ and 
	\begin{equation*}
		\left\|\boldsymbol{Z}_{\scriptscriptstyle \xi_i^{(k)}}\right\|_2 \leq \frac{\max\limits_{j=1,\ldots,n}\left\{\frac{\|\boldsymbol{a}_j\|_{2}^{2}}{p_j}\right\}-\lambda_d\left(\boldsymbol{A}^{\top} \boldsymbol{A}\right)}{B},~i=1,\ldots,B.
	\end{equation*}
	
	The next step is to estimate the matrix variance statistic \eqref{eq: matrix variance statistic} of $\boldsymbol{Z}_k$. Firstly, we have
	\begin{equation*}
		\textstyle
		\boldsymbol{Z}_{\scriptscriptstyle \xi_i^{(k)}}^{\top} \boldsymbol{Z}_{\scriptscriptstyle \xi_i^{(k)}}=\frac{1}{B^2}\left(\frac{\left\|\boldsymbol{a}_{\scriptscriptstyle \xi_i^{(k)}}\right\|_2^2}{p_{\scriptscriptstyle \xi_i^{(k)}}^2} \boldsymbol{a}_{\scriptscriptstyle \xi_i^{(k)}} \boldsymbol{a}_{\scriptscriptstyle \xi_i^{(k)}}^{\top}-\frac{1}{p_{\scriptscriptstyle \xi_i^{(k)}}\!} \boldsymbol{A}^{\top} \boldsymbol{A} \boldsymbol{a}_{\scriptscriptstyle \xi_i^{(k)}} \boldsymbol{a}_{\scriptscriptstyle \xi_i^{(k)}}^{\top}-\frac{1}{p_{\scriptscriptstyle \xi_i^{(k)}}\!} \boldsymbol{a}_{\scriptscriptstyle \xi_i^{(k)}} \boldsymbol{a}_{\scriptscriptstyle \xi_i^{(k)}}^{\top} \boldsymbol{A}^{\top} \boldsymbol{A}+\left(\boldsymbol{A}^{\top} \boldsymbol{A}\right)^2 \right)
	\end{equation*}
	for each $i=1,\ldots,B$. Then, the expectation of $\boldsymbol{Z}_{\scriptscriptstyle \xi_i^{(k)}}^{\top} \boldsymbol{Z}_{\scriptscriptstyle \xi_i^{(k)}} $ can be estimated as
	\begin{equation*}
		\mathbb{E}\left[\boldsymbol{Z}_{\scriptscriptstyle \xi_i^{(k)}}^{\top} \boldsymbol{Z}_{\scriptscriptstyle \xi_i^{(k)}} \!\right]\!=\!\frac{1}{B^2}\!\left( \sum_{i=1}^n \frac{\left\|\boldsymbol{a}_i\right\|_2^2}{p_i} \boldsymbol{a}_i \boldsymbol{a}_i^{\top}-\left(\boldsymbol{A}^{\top} \!\boldsymbol{A}\right)^2\right)\!\preceq\! \frac{1}{B^2}\!\left(\max\limits_{j=1,\ldots,n}\left\{\frac{\|\boldsymbol{a}_j\|_{2}^{2}}{p_j}\right\} \boldsymbol{A}^{\top}\! \boldsymbol{A}-\left(\boldsymbol{A}^{\top}\! \boldsymbol{A}\right)^2\right)\!,
	\end{equation*}
	where we use $\boldsymbol{A}^{\top}\boldsymbol{A}=\sum_{i=1}^n\boldsymbol{a}_i\boldsymbol{a}_i^{\top}$ in the first equality. Notice that $\boldsymbol{Z}_{\scriptscriptstyle \xi_i^{(k)}}^{\top} \boldsymbol{Z}_{\scriptscriptstyle \xi_i^{(k)}}$ is positive semidefinite. The above implies that for each $i=1\ldots,B$,
	\begin{align*}
		\left\|\mathbb{E}\left[\boldsymbol{Z}_{\scriptscriptstyle \xi_i^{(k)}}^{\top} \boldsymbol{Z}_{\scriptscriptstyle \xi_i^{(k)}}\right]\right\|_2 &\! \leq\! \frac{\left\|\max\limits_{j=1,\ldots,n}\left\{\frac{\|\boldsymbol{a}_j\|_{2}^{2}}{p_j}\right\} \boldsymbol{A}^{\top}\! \boldsymbol{A}\!-\!\left(\boldsymbol{A}^{\top} \!\boldsymbol{A}\right)^2\right\|_2}{B^2} \nonumber \\
		& \!\leq\! \frac{\left\|\boldsymbol{A}^{\top} \!\boldsymbol{A}\right\|_2\left\|\max\limits_{j=1,\ldots,n}\!\left\{\!\frac{\|\boldsymbol{a}_j\|_{2}^{2}}{p_j}\!\right\} \!\boldsymbol{I}\!-\!\boldsymbol{A}^{\top}\! \boldsymbol{A}\right\|_2}{B^2}\!=\! \frac{\|\boldsymbol{A}\|_2^2\!\left(\!\max\limits_{j=1,\ldots,n}\!\left\{\!\frac{\|\boldsymbol{a}_j\|_{2}^{2}}{p_j}\!\right\}\!-\!\lambda_d\!\left(\boldsymbol{A}^{\top}\! \boldsymbol{A}\right)\hspace{-0.4em}\right)}{B^2},
	\end{align*}
	where the last equality is derived from \cref{lem:probability bound}, which implies that $\max\limits_{j=1,\ldots,n}\left\{\|\boldsymbol{a}_j\|_{2}^{2}/p_j\right\} \boldsymbol{I}-\boldsymbol{A}^{\top} \boldsymbol{A}$ is positive semidefinite. As a result, the matrix variance statistic of $\boldsymbol{Z}_k$ satisfies
	\begin{equation*}
		\nu(\boldsymbol{Z}_k)\!=\!\left\|\sum_{i=1}^B \mathbb{E}\left[\boldsymbol{Z}_{\scriptscriptstyle \xi_i^{(k)}}^{\top} \boldsymbol{Z}_{\scriptscriptstyle \xi_i^{(k)}}\right]\right\|_2\! \leq \!\sum_{i=1}^{B} \left\| \mathbb{E}\left[\boldsymbol{Z}_{\scriptscriptstyle \xi_i^{(k)}}^{\top} \boldsymbol{Z}_{\scriptscriptstyle \xi_i^{(k)}}\right]\right\|_2 \!\leq\! \frac{\|\boldsymbol{A}\|_2^2\left(\!\max\limits_{j=1,\ldots,n}\!\left\{\!\frac{\|\boldsymbol{a}_j\|_{2}^{2}}{p_j}\!\right\}-\lambda_d(\boldsymbol{A}^{\top} \boldsymbol{A})\!\right)}{B}.
	\end{equation*}
	
	Combined with the independence of the matrices $\{\boldsymbol{Z}_{\scriptscriptstyle \xi_i^{(k)}}\}_{i=1}^{B}$ and $\mathbb{E}[\boldsymbol{Z}_{\scriptscriptstyle \xi_i^{(k)}}]=\boldsymbol{O}$, applying \cref{lem:lemma Matrix Bernstein}, we obtain \eqref{eq: standard error upper bound}.
\end{proof}

The next theorem shows an R-linear convergence rate of \hyperref[alg:SGA-RMSProp]{SGA-RMSProp} on the consistent LLSP, which is our main result.
\begin{theorem} \label{thm: convergence of consistent LLSP}
	Let $\boldsymbol{x}^*$ be the optimal point of the consistent LLSP \eqref{eq:LLSP obj}, and $\{\boldsymbol{x}_k\}_{k=1}^{K}$ be the sequence generated by \hyperref[alg:SGA-RMSProp]{SGA-RMSProp} with $\eta_k\equiv \eta:=\frac{2}{\overline{u}\left(\lambda_1\left(\boldsymbol{A}^{\top}\boldsymbol{A}\right)+\lambda_d\left(\boldsymbol{A}^{\top}\boldsymbol{A}\right)\right)}$ . Then we have
	
	\begin{equation}\label{eq: linear convergence rate initial}
		\mathbb{E}\left[\left\|\boldsymbol{x}_{K+1}-\boldsymbol{x}^*\right\|_2\right]\leq \rho\left\|\boldsymbol{x}_1-\boldsymbol{x}^*\right\|_2 \gamma^K,
	\end{equation}
	where $\rho\leq1$ and
	\begin{equation}\label{eq: gamma}
		\gamma:=\text{G}(\gamma_1,\ldots,\gamma_K)\left(1-\frac{2 \left(\underline{u} \lambda_d\left(\boldsymbol{A}^{\top}\boldsymbol{A}\right)-\overline{u} \sigma\right)}{ \overline{u}\left(\lambda_1\left(\boldsymbol{A}^{\top}\boldsymbol{A}\right)+\lambda_d\left(\boldsymbol{A}^{\top}\boldsymbol{A}\right)\right)}+\varepsilon\right),
	\end{equation}
	in which $\gamma_k\leq1$, $k=1,\ldots ,K$, $\text{G}(\gamma_1,\ldots,\gamma_K):=\big(\prod\limits_{k=1}^{K}\gamma_k\big)^{\frac{1}{K}}$, and
	\begin{align}\label{eq: sigma}
		\textstyle
		\sigma :=& \Bigg(\frac{2\left(\max\limits_{j=1,\ldots,n}\left\{\frac{\|\boldsymbol{a}_j\|_{2}^{2}}{p_j}\right\}-\lambda_d\left(\boldsymbol{A}^{\top} \boldsymbol{A}\right)\right)\|\boldsymbol{A}\|_2^2 \log (2 d)}{B}\Bigg)^{\frac{1}{2}}
		\nonumber \\
		&+ \frac{\log (2 d)}{3B}\left(\max\limits_{\scriptscriptstyle j=1,\ldots,n}\left\{\frac{\|\boldsymbol{a}_j\|_{2}^{2}}{p_j}\right\}-\lambda_d\left(\boldsymbol{A}^{\top} \boldsymbol{A}\right)\right).
	\end{align}
\end{theorem}

\begin{proof}
	Revisiting \eqref{eq: norm inequality}, noticing that $\boldsymbol{D}_k$ is a diagonal matrix for $k=0,\ldots,K$, which makes their multiplication commuted, we have
	\begin{align}
		& \left\|\boldsymbol{Y}_K  \cdots \boldsymbol{Y}_1\right\|_2 \nonumber \\
		= & \left\|\left(\boldsymbol{I}-\eta \boldsymbol{D}_K \boldsymbol{M}_K\right)\left(\boldsymbol{I}-\eta \boldsymbol{D}_{K-1} \boldsymbol{M}_{K-1}\right) \cdots\left(\boldsymbol{I}-\eta \boldsymbol{D}_1 \boldsymbol{M}_1\right)\right\|_2 \nonumber \\
		= & \left\|\boldsymbol{D}_K\left(\boldsymbol{D}_K^{-1} \boldsymbol{D}_{K-1}-\eta \boldsymbol{M}_K \boldsymbol{D}_{K-1}\right) \boldsymbol{D}_{K-1}^{-1} \boldsymbol{D}_{K-1}\left(\boldsymbol{D}_{K-1}^{-1} \boldsymbol{D}_{K-2}-\eta \boldsymbol{M}_{K-1} \boldsymbol{D}_{K-2}\right) \right. \nonumber \\
		& \left. \cdots\left(\boldsymbol{D}_1^{-1} \boldsymbol{D}_0-\eta \boldsymbol{M}_1 \boldsymbol{D}_0\right) \boldsymbol{D}_0^{-1}\right\|_2 \nonumber \\
		= & {\small\left\|\boldsymbol{D}_K \boldsymbol{D}_{K-1}^{-\frac{1}{2}}\left(\boldsymbol{D}_K^{-1} \boldsymbol{D}_{K-1}-\eta \boldsymbol{D}_{K-1}^{\frac{1}{2}} \boldsymbol{M}_K \boldsymbol{D}_{K-1}^{\frac{1}{2}}\right) \boldsymbol{D}_{K-1}^{\frac{1}{2}} \boldsymbol{D}_{K-2}^{-\frac{1}{2}} \right.} \nonumber \\
		& {\small \left. \left(\boldsymbol{D}_{K-1}^{-1} \boldsymbol{D}_{K-2}-\eta \boldsymbol{D}_{K-2}^{\frac{1}{2}} \boldsymbol{M}_{K-1} \boldsymbol{D}_{K-2}^{\frac{1}{2}}\right) \boldsymbol{D}_{K-2}^{\frac{1}{2}} \cdots \boldsymbol{D}_0^{-\frac{1}{2}}\left(\boldsymbol{D}_1^{-1} \boldsymbol{D}_0-\eta \boldsymbol{D}_0^{\frac{1}{2}} \boldsymbol{M}_1 \boldsymbol{D}_0^{\frac{1}{2}}\right) \boldsymbol{D}_0^{-\frac{1}{2}}\right\|_2} \nonumber \\
		\leq &{\small \left\|\boldsymbol{D}_K^{\frac{1}{2}}\right\|_2\left\|\boldsymbol{D}_K^{\frac{1}{2}} \boldsymbol{D}_{K-1}^{-\frac{1}{2}}\right\|_2\left\|\boldsymbol{D}_K^{-1} \boldsymbol{D}_{K-1}-\eta \boldsymbol{D}_{K-1}^{\frac{1}{2}} \boldsymbol{M}_K \boldsymbol{D}_{K-1}^{\frac{1}{2}}\right\|_2 \left\|\boldsymbol{D}_{K-1}^{\frac{1}{2}} \boldsymbol{D}_{K-2}^{-\frac{1}{2}}\right\|_2}\nonumber \\
		& {\small\left\|\boldsymbol{D}_{\!K\!-\!1}^{\!-\!1} \!\boldsymbol{D}_{\!K\!-\!2}\!-\!\eta \boldsymbol{D}_{\!K\!-\!2}^{\!\frac{1}{2}} \boldsymbol{M}_{\!K\!-\!1}\! \boldsymbol{D}_{\!K\!-\!2}^{\!\frac{1}{2}}\right\|_2\! \left\|\boldsymbol{D}_{\!K\!-\!2}^{\!\frac{1}{2}}\!\boldsymbol{D}_{\!K\!-\!3}^{\!-\!\frac{1}{2}}\right\|_2  \!\!\!\cdots\!\left\|\boldsymbol{D}_{\!1}^{\!\frac{1}{2}}\! \boldsymbol{D}_{\!0}^{\!-\!\frac{1}{2}}\right\|_2\!\left\|\boldsymbol{D}_{\!1}^{\!-\!1} \!\boldsymbol{D}_{\!0}\!-\!\eta \boldsymbol{D}_{\!0}^{\!\frac{1}{2}} \!\boldsymbol{M}_{\!1}\! \boldsymbol{D}_{\!0}^{\!\frac{1}{2}}\right\|_2\!\left\|\boldsymbol{D}_{\!0}^{\!-\!\frac{1}{2}}\right\|_2.}\label{eq: Y transfer}
	\end{align}
	
	Denote $\boldsymbol{X}_k:=\boldsymbol{I}-\eta \boldsymbol{D}_{k-1}^{\frac{1}{2}} \boldsymbol{M}_k \boldsymbol{D}_{k-1}^{\frac{1}{2}}$, $k=1,\ldots,K$. \cref{lem:SGA-RMSProp control} implies that
	\begin{equation}\label{eq: D_k to X_k}
		\left\|\boldsymbol{D}_{\!k}^{-\!1} \boldsymbol{D}_{\!k\!-\!1}\!-\!\eta \boldsymbol{D}_{\!k\!-\!1}^{\!\frac{1}{2}} \boldsymbol{M}_{\!k} \boldsymbol{D}_{\!k\!-\!1}^{\!\frac{1}{2}}\right\|_2\!\!\! \leq\!\left\|\boldsymbol{I}\!-\!\eta \boldsymbol{D}_{\!k\!-\!1}^{\!\frac{1}{2}} \boldsymbol{M}_{\!k} \boldsymbol{D}_{\!k\!-\!1}^{\!\frac{1}{2}}\right\|_2\!\!\!+\!\left\|\boldsymbol{D}_{\!k}^{-\!1} \boldsymbol{D}_{\!k\!-\!1}\!-\!\boldsymbol{I}\right\|_2  \!\!\overset{\eqref{eq:D_k D_k-1 estimate}}{\leq}\!\!\left\|\boldsymbol{X}_{\!k}\right\|_2\!+\!\varepsilon.
	\end{equation}
	
	According to the setting $\|\boldsymbol{D}_0^{-\frac{1}{2}}\|_2=1/\sqrt{\overline{u}}$ and \cref{lem:SGA-RMSProp control}, we have that $\|\boldsymbol{D}_{\!K}^{\!\frac{1}{2}}\!\|_2 \|\boldsymbol{D}_{\!0}^{\!-\!\frac{1}{2}}\!\|_2$ and $\|\boldsymbol{D}_k^{\frac{1}{2}}\boldsymbol{D}_{k-1}^{-\frac{1}{2}}\|_2$, $k=1,\ldots,K$, are bounded above by some constants at most 1. These constants, denoted by $\rho$, $\gamma_k$, $k=1,\ldots,K$ respectively, only depend on the hyperparameters of \hyperref[alg:SGA-RMSProp]{SGA-RMSProp}, the initial point $\boldsymbol{x}_1$, and $\boldsymbol{A}$, $\boldsymbol{b}$ in LLSP. Recalling that $\mathcal{F}_k$ is the $\sigma$-field generated by the random variables in $\mathcal{S}_1,\ldots,\mathcal{S}_{k-1}$, $k=2,\ldots,K$, taking the expectation for both sides of the inequality \eqref{eq: Y transfer} and by \cref{prop:CE tower}, we have 
	{\small
	\begin{align}
		\mathbb{E}\left[\left\|\boldsymbol{Y}_K \cdots \boldsymbol{Y}_1\right\|_2\right] &\!\! \overset{\eqref{eq: Y transfer}}{\leq}\!\! \rho\! \left(\prod_{k=1}^{K}\gamma_k\right) \mathbb{E}\left[\prod_{k=1}^{K}\left\|\boldsymbol{D}_k^{-1} \boldsymbol{D}_{k-1}-\eta \boldsymbol{D}_{k-1}^{\frac{1}{2}} \boldsymbol{M}_k \boldsymbol{D}_{k-1}^{\frac{1}{2}}\right\|_2\right] \nonumber \\
		& \!\!\overset{\eqref{eq: D_k to X_k}}{\leq}\!\! \rho\! \left(\prod_{k=1}^{K}\!\!\gamma_k\!\!\right) \!\mathbb{E}\!\left[\prod_{k=1}^K\!\!\left(\left\|\boldsymbol{X}_{\!k}\right\|_2\!\!+\!\varepsilon\right)\right]\!\! \overset{\eqref{eq: CE tower}}{=}\!\!\rho \!\left(\prod_{k=1}^{K}\!\!\gamma_k\!\!\right)\! \mathbb{E}\!\left[\!\cdots\!\mathbb{E}\!\left[\mathbb{E}\!\left[\prod_{k=1}^K\!\!\left(\left\|\boldsymbol{X}_{\!k}\right\|_2\!\!+\!\varepsilon\right) \!\Bigg|\mathcal{F}_{\!K}\! \right]\! \Bigg| \mathcal{F}_{\!K\!-\!1}\right]\! \!\cdots\! \right]. \label{eq: law of total expectation}
	\end{align}
}
	
	By $\boldsymbol{D}_{k-1}=\text{diag}(\frac{1}{\sqrt{u_{k-1,1}}},\ldots,\frac{1}{\sqrt{u_{k-1,d}}})$ and \eqref{eq: RMSProp u iter}, $\boldsymbol{D}_{k-1}$ is $\mathcal{F}_K$-measurable for $k=1,\ldots,K$. Moreover, the matrix $\boldsymbol{M}_k$ depends solely on $\mathcal{S}_k$ and $\boldsymbol{X}_k=\boldsymbol{I}-\eta \boldsymbol{D}_{k-1}^{\frac{1}{2}} \boldsymbol{M}_k \boldsymbol{D}_{k-1}^{\frac{1}{2}}$, implying that $\boldsymbol{X}_k$ is $\mathcal{F}_K$-measurable for $k=1,\ldots,K-1$. Therefore, we have
	\begin{equation}\label{eq: measurable expectation formula}
		\mathbb{E}\left[\prod_{k=1}^K\left(\left\|\boldsymbol{X}_k\right\|_2+\varepsilon\right)\Bigg | \mathcal{F}_K\right]\overset{\eqref{eq: CE measurable}}{=}\mathbb{E}\left[\left\|\boldsymbol{X}_K\right\|_2+\varepsilon\big|\mathcal{F}_K\right]\prod_{k=1}^{K-1}\left(\left\|\boldsymbol{X}_k\right\|_2+\varepsilon\right).
	\end{equation}
	Using triangle inequality, $\mathbb{E}\left[\left\|\boldsymbol{X}_K\right\|_2+\varepsilon|\mathcal{F}_K\right]$ has the following upper bound estimation
	\begin{align}\label{eq: expectation triangle inequality}
		\mathbb{E}\left[\left\|\boldsymbol{X}_K\right\|_2+\varepsilon|\mathcal{F}_K\right]=&\mathbb{E}\left[\left\|\boldsymbol{X}_K-\mathbb{E}\left[\boldsymbol{X}_K|\mathcal{F}_K\right]+\mathbb{E}\left[\boldsymbol{X}_K|\mathcal{F}_K\right]\right\|_2|\mathcal{F}_K\right]+\varepsilon \nonumber \\
		\leq & \mathbb{E} \left[\left\|\boldsymbol{X}_K-\mathbb{E}\left[\boldsymbol{X}_K|\mathcal{F}_K\right]\right\|_2+\left\|\mathbb{E}\left[\boldsymbol{X}_K|\mathcal{F}_K\right]\right\|_2|\mathcal{F}_K\right] +\varepsilon \nonumber\\
		=&\mathbb{E}\left[\left\|\boldsymbol{X}_K-\mathbb{E}\left[\boldsymbol{X}_K|\mathcal{F}_K\right]\right\|_2|\mathcal{F}_K\right]+\left\|\mathbb{E}\left[\boldsymbol{X}_K|\mathcal{F}_K\right]\right\|_2+\varepsilon.
	\end{align}
	
	As $\boldsymbol{D}_{K-1}$ is $\mathcal{F}_K$-measurable and $\boldsymbol{M}_K$ is independent of $\mathcal{F}_K$, we have
	\begin{align}\label{eq: calculate conditional expectation in norm}
		\left\|\mathbb{E}\left[\boldsymbol{X}_K|\mathcal{F}_K\right]\right\|_2&\hspace{5mm}=\hspace{5mm}\left\|\mathbb{E}\left[\boldsymbol{I}-\eta \boldsymbol{D}_{K-1}^{\frac{1}{2}} \boldsymbol{M}_K \boldsymbol{D}_{K-1}^{\frac{1}{2}}\bigg|\mathcal{F}_K\right]\right\|_2 \nonumber \\
		&\overset{\eqref{eq: CE measurable},\eqref{eq: CE independent}}{=}\left\|\boldsymbol{I}-\eta \boldsymbol{D}_{K-1}^{\frac{1}{2}} \mathbb{E}\left[ \boldsymbol{M}_K\right] \boldsymbol{D}_{K-1}^{\frac{1}{2}}\right\|_2 =\left\|\boldsymbol{I}-\eta \boldsymbol{D}_{K-1}^{\frac{1}{2}} \boldsymbol{A}^{\top} \boldsymbol{A} \boldsymbol{D}_{K-1}^{\frac{1}{2}}\right\|_2.
	\end{align}
	Since $\boldsymbol{I}-\eta \boldsymbol{D}_{K-1}^{\frac{1}{2}} \boldsymbol{A}^{\top} \boldsymbol{A} \boldsymbol{D}_{K-1}^{\frac{1}{2}}$ is a
	symmetric matrix, the spectral norm can be calculated as
	\begin{align}\label{eq: calculate spectral norm}
		\left\|\boldsymbol{I}\!-\!\eta \boldsymbol{D}_{K-1}^{\frac{1}{2}} \boldsymbol{A}^{\top}\! \boldsymbol{A} \boldsymbol{D}_{K-1}^{\frac{1}{2}}\right\|_2 \!\!=&\max\! \left\{\!\lambda_1\!\left(\!\boldsymbol{I}\!-\!\eta \boldsymbol{D}_{K-1}^{\frac{1}{2}} \boldsymbol{A}^{\top} \boldsymbol{A} \boldsymbol{D}_{K-1}^{\frac{1}{2}}\!\right),-\lambda_d\!\left(\!\boldsymbol{I}\!-\!\eta \boldsymbol{D}_{K-1}^{\frac{1}{2}} \boldsymbol{A}^{\top}\! \boldsymbol{A} \boldsymbol{D}_{K-1}^{\frac{1}{2}}\!\right)\!\right\} \nonumber \\
		\!\!=&\max\! \left\{\!1\!-\!\eta \lambda_d\!\left(\!\boldsymbol{D}_{K-1}^{\frac{1}{2}} \boldsymbol{A}^{\top}\! \boldsymbol{A} \boldsymbol{D}_{K-1}^{\frac{1}{2}}\!\right),\eta \lambda_1\!\left(\! \boldsymbol{D}_{K-1}^{\frac{1}{2}} \boldsymbol{A}^{\top}\! \boldsymbol{A} \boldsymbol{D}_{K-1}^{\frac{1}{2}}\!\right)\!-\!1\!\right\}.
	\end{align}
	Note that when a matrix is positive semidefinite, its singular values are the same as its eigenvalues. By \cref{lem:Ky Fan} and \cref{lem:SGA-RMSProp control}, we have
	{\small
	\begin{equation}\label{eq: lambda_1 upper bound}
		\lambda_1\!\left(\!\boldsymbol{D}_{\!K\!-\!1}^{\!\frac{1}{2}} \boldsymbol{A}^{\!\top}\! \boldsymbol{A} \boldsymbol{D}_{\!K\!-\!1}^{\!\frac{1}{2}}\!\right)\!\!=\!\!s_1\!\left(\!\boldsymbol{D}_{\!K\!-\!1}^{\!\frac{1}{2}} \boldsymbol{A}^{\!\top}\! \boldsymbol{A} \boldsymbol{D}_{\!K\!-\!1}^{\!\frac{1}{2}}\!\right)\! \!\overset{\eqref{lem: Fan inequality}}{\leq}\!\! s_1\!\left(\!\boldsymbol{D}_{\!K\!-\!1}^{\!\frac{1}{2}}\!\right)\!s_1\!\left(\!\boldsymbol{A}^{\!\top}\! \boldsymbol{A}\!\right)s_1\!\left(\!\boldsymbol{D}_{\!K\!-\!1}^{\!\frac{1}{2}}\!\right)\!\!\overset{\eqref{eq: D_k estimate}}{\leq}\!\! \overline{u}\lambda_1\!\left(\!\boldsymbol{A}^{\!\top}\! \boldsymbol{A}\!\right),
	\end{equation}
}
	and
	{\small
	\begin{equation}\label{eq: lambda_d upper bound}
		\lambda_d\!\left(\!\boldsymbol{D}_{\!K\!-\!1}^{\!\frac{1}{2}} \boldsymbol{A}^{\!\top}\! \boldsymbol{A} \boldsymbol{D}_{\!K\!-\!1}^{\!\frac{1}{2}}\!\right)\!\!=\!\!s_d\!\left(\!\boldsymbol{D}_{\!K\!-\!1}^{\!\frac{1}{2}} \boldsymbol{A}^{\!\top}\! \boldsymbol{A} \boldsymbol{D}_{\!K\!-\!1}^{\!\frac{1}{2}}\!\right) \!\!\overset{\eqref{lem: Fan inequality}}{\leq}\!\! s_1\!\left(\!\boldsymbol{D}_{\!K\!-\!1}^{\!\frac{1}{2}}\!\right)s_d\!\left(\!\boldsymbol{A}^{\!\top} \!\boldsymbol{A}\!\right)s_1\!\left(\!\boldsymbol{D}_{\!K\!-\!1}^{\!\frac{1}{2}}\!\right)\!\!\overset{\eqref{eq: D_k estimate}}{\leq}\!\! \overline{u}\lambda_d\!\left(\!\boldsymbol{A}^{\!\top} \!\boldsymbol{A}\!\right).
	\end{equation}
}
	Combined with our step size setting $\eta=\frac{2}{\overline{u}\left(\lambda_1\left(\boldsymbol{A}^{\top}\boldsymbol{A}\right)+\lambda_d\left(\boldsymbol{A}^{\top}\boldsymbol{A}\right)\right)}$, we have
	\begin{align}
		1-\eta \lambda_d\left(\boldsymbol{D}_{K-1}^{\frac{1}{2}} \boldsymbol{A}^{\top} \boldsymbol{A} \boldsymbol{D}_{K-1}^{\frac{1}{2}}\right)&\hspace{2.3mm}=\hspace{2.3mm}\frac{\overline{u}\left(\lambda_1\left(\boldsymbol{A}^{\top}\boldsymbol{A}\right)+\lambda_d\left(\boldsymbol{A}^{\top}\boldsymbol{A}\right)\right)-2\lambda_d\left(\boldsymbol{D}_{K-1}^{\frac{1}{2}} \boldsymbol{A}^{\top} \boldsymbol{A} \boldsymbol{D}_{K-1}^{\frac{1}{2}}\right)}{\overline{u}\left(\lambda_1\left(\boldsymbol{A}^{\top}\boldsymbol{A}\right)+\lambda_d\left(\boldsymbol{A}^{\top}\boldsymbol{A}\right)\right)} \nonumber \\
		&\overset{\eqref{eq: lambda_d upper bound}}{\geq}\frac{\overline{u}\left(\lambda_1\left(\boldsymbol{A}^{\top}\boldsymbol{A}\right)+\lambda_d\left(\boldsymbol{A}^{\top}\boldsymbol{A}\right)\right)-2\overline{u}\lambda_d\left(\boldsymbol{A}^{\top} \boldsymbol{A}\right)}{\overline{u}\left(\lambda_1\left(\boldsymbol{A}^{\top}\boldsymbol{A}\right)+\lambda_d\left(\boldsymbol{A}^{\top}\boldsymbol{A}\right)\right)} \nonumber \\
		&\hspace{2.3mm}=\hspace{2.3mm}\frac{2 \overline{u}\lambda_1\left(\boldsymbol{A}^{\top}\boldsymbol{A}\right)-\overline{u}\left(\lambda_1\left(\boldsymbol{A}^{\top}\boldsymbol{A}\right)+\lambda_d\left(\boldsymbol{A}^{\top}\boldsymbol{A}\right)\right)}{\overline{u}\left(\lambda_1\left(\boldsymbol{A}^{\top}\boldsymbol{A}\right)+\lambda_d\left(\boldsymbol{A}^{\top}\boldsymbol{A}\right)\right)} \nonumber \\
		&\overset{\eqref{eq: lambda_1 upper bound}}{\geq}\frac{2 \lambda_1\left(\boldsymbol{D}_{K-1}^{\frac{1}{2}} \boldsymbol{A}^{\top} \boldsymbol{A} \boldsymbol{D}_{K-1}^{\frac{1}{2}}\right)-\overline{u}\left(\lambda_1\left(\boldsymbol{A}^{\top}\boldsymbol{A}\right)+\lambda_d\left(\boldsymbol{A}^{\top}\boldsymbol{A}\right)\right)}{\overline{u}\left(\lambda_1\left(\boldsymbol{A}^{\top}\boldsymbol{A}\right)+\lambda_d\left(\boldsymbol{A}^{\top}\boldsymbol{A}\right)\right)} \nonumber \\
		&\hspace{2.3mm}=\hspace{2.3mm}\eta \lambda_1\left( \boldsymbol{D}_{K-1}^{\frac{1}{2}} \boldsymbol{A}^{\top} \boldsymbol{A} \boldsymbol{D}_{K-1}^{\frac{1}{2}}\right)-1. \nonumber
	\end{align} 	
	Thus, \eqref{eq: calculate spectral norm} implies that
	\begin{equation}\label{eq: calculate spectral norm final}
		\left\|\boldsymbol{I}-\eta \boldsymbol{D}_{K-1}^{\frac{1}{2}} \boldsymbol{A}^{\top} \boldsymbol{A} \boldsymbol{D}_{K-1}^{\frac{1}{2}}\right\|_2 =1-\eta \lambda_d\left(\boldsymbol{D}_{K-1}^{\frac{1}{2}} \boldsymbol{A}^{\top} \boldsymbol{A} \boldsymbol{D}_{K-1}^{\frac{1}{2}}\right).
	\end{equation}
	\cref{lem:Ky Fan} and \cref{lem:SGA-RMSProp control} also imply that
	\begin{align}\label{eq: upper bound of lambda_d AA}
		\lambda_d\left(\boldsymbol{A}^{\top} \boldsymbol{A}\right) & \hspace{1.4mm}=\hspace{1.4mm}\lambda_d\left(\boldsymbol{D}_{K-1}^{-\frac{1}{2}} \boldsymbol{D}_{K-1}^{\frac{1}{2}} \boldsymbol{A}^{\top} \boldsymbol{A} \boldsymbol{D}_{K-1}^{\frac{1}{2}} \boldsymbol{D}_{K-1}^{-\frac{1}{2}}\right) \nonumber \\
		&\small \overset{\eqref{lem: Fan inequality}}{\leq} s_1\left(\boldsymbol{D}_{K-1}^{-\frac{1}{2}}\right) s_d\left(\boldsymbol{D}_{K-1}^{\frac{1}{2}} \boldsymbol{A}^{\top} \boldsymbol{A} \boldsymbol{D}_{K-1}^{\frac{1}{2}} \right)s_1\left(\boldsymbol{D}_{K-1}^{-\frac{1}{2}}\right)\overset{\eqref{eq: D_k estimate}}{\leq} \frac{\lambda_d\left(\boldsymbol{D}_{K-1}^{\frac{1}{2}} \boldsymbol{A}^{\top} \boldsymbol{A} \boldsymbol{D}_{K-1}^{\frac{1}{2}}\right)}{\underline{u}} .
	\end{align}
	Combining \eqref{eq: calculate conditional expectation in norm}, \eqref{eq: calculate spectral norm final} with \eqref{eq: upper bound of lambda_d AA}, we have
	\begin{equation}\label{eq: expectation upper bound}
		\left\|\mathbb{E}\left[\boldsymbol{X}_K|\mathcal{F}_K\right]\right\|_2 \leq 1-\eta \underline{u} \lambda_d\left(\boldsymbol{A}^{\top} \boldsymbol{A}\right)=1- \frac{2\underline{u} \lambda_d\left(\boldsymbol{A}^{\top}\boldsymbol{A}\right)}{ \overline{u}\left(\lambda_1\left(\boldsymbol{A}^{\top}\boldsymbol{A}\right)+\lambda_d\left(\boldsymbol{A}^{\top}\boldsymbol{A}\right)\right)}.
	\end{equation}
	
	Moreover,  given $\sigma$ as in \eqref{eq: sigma}, by \cref{lem:SGA-RMSProp control} and \cref{lem: standard error upper bound}, we have
	{\small
	\begin{align}\label{eq: second part}
			&\mathbb{E}\left[\left\|\boldsymbol{X}_K-\mathbb{E}\left[\boldsymbol{X}_K|\mathcal{F}_K\right]\right\|_2|\mathcal{F}_K\right] \nonumber \\
			= \!\hspace{1.6mm}& \eta \mathbb{E}\left[\left\|\boldsymbol{D}_{K-1}^{\frac{1}{2}}\left(-\boldsymbol{M}_K+\boldsymbol{A}^{\top} \boldsymbol{A}\right) \boldsymbol{D}_{K-1}^{\frac{1}{2}}\right\|_2\bigg|\mathcal{F}_K\right] \nonumber \\
			\overset{\eqref{eq: D_k estimate}}{\leq}\!& \eta \overline{u} \mathbb{E}\left[\left\|\left(-\boldsymbol{M}_K+\boldsymbol{A}^{\top} \boldsymbol{A}\right)\right\|_2\Big|\mathcal{F}_K\right] \nonumber \\
			\overset{\eqref{eq: CE independent}}{=} \!& \eta \overline{u} \mathbb{E}\left[\left\|\sum_{i=1}^B \frac{1}{B}\left(-\frac{1}{p_{\scriptscriptstyle \xi_{i}^{(K)}}} \boldsymbol{a}_{\scriptscriptstyle \xi_i^{(K)}} \boldsymbol{a}_{\scriptscriptstyle \xi_i^{(K)}}^{\top}+\boldsymbol{A}^{\top} \boldsymbol{A}\right)\right\|_2\right] \nonumber \\
			\overset{\eqref{eq: standard error upper bound}}{\leq}\!&
				 \eta \overline{u} \Bigg(\!\!\Bigg(\!\frac{2\Big(\!\max\limits_{\scriptscriptstyle j=1,\ldots,n}\!\left\{\!\frac{\|\boldsymbol{a}_j\|_{2}^{2}}{p_j}\!\right\}\!-\!\lambda_d\!\left(\!\boldsymbol{A}^{\!\top}\! \boldsymbol{A}\!\right)\!\!\Big)\|\boldsymbol{A}\|_2^2 \log (2 d)}{B}\!\Bigg)^{\!\!\frac{1}{2}}\!\!+\!\!\frac{\log(2d)\Big(\!\max\limits_{\scriptscriptstyle j=1,\ldots,n}\!\left\{\!\frac{\|\boldsymbol{a}_j\|_{2}^{2}}{p_j}\!\right\}\!-\!\lambda_d\left(\!\boldsymbol{A}^{\!\top}\! \boldsymbol{A}\!\right)\!\!\Big)}{3B}\!\Bigg)\nonumber \\
			 \overset{\eqref{eq: sigma}}{=} \!&   \frac{2\sigma}{\lambda_1\left(\boldsymbol{A}^{\top} \boldsymbol{A}\right)+\lambda_d\left(\boldsymbol{A}^{\top} \boldsymbol{A}\right)}.
		\end{align}
	}
	Substituting \eqref{eq: second part} and \eqref{eq: expectation upper bound} into \eqref{eq: expectation triangle inequality}, the estimation of $\mathbb{E}[\prod_{k=1}^K(\left\|\boldsymbol{X}_k\right\|_2+\varepsilon)|\mathcal{F}_K]$ by \eqref{eq: measurable expectation formula} becomes
	\begin{equation}\label{eq: X_k estimate}
		\textstyle
		\mathbb{E}\left[\prod\limits_{k=1}^K\left(\left\|\boldsymbol{X}_k\right\|_2+\varepsilon\right)\bigg|\mathcal{F}_K\right] \leq\left(1- \frac{2 \left(\underline{u} \lambda_d\left(\boldsymbol{A}^{\top}\boldsymbol{A}\right)-\overline{u} \sigma\right)}{ \overline{u}\left(\lambda_1\left(\boldsymbol{A}^{\top}\boldsymbol{A}\right)+\lambda_d\left(\boldsymbol{A}^{\top}\boldsymbol{A}\right)\right)}+\varepsilon\right) \prod\limits_{k=1}^{K-1}\left(\left\|\boldsymbol{X}_k\right\|_2+\varepsilon\right).
	\end{equation}
	
	Notice that $1- \frac{2 \left(\underline{u} \lambda_d\left(\boldsymbol{A}^{\top}\boldsymbol{A}\right)-\overline{u} \sigma\right)}{ \overline{u}\left(\lambda_1\left(\boldsymbol{A}^{\top}\boldsymbol{A}\right)+\lambda_d\left(\boldsymbol{A}^{\top}\boldsymbol{A}\right)\right)}+\varepsilon$ has already been a constant, so we can extend the same procedure to $k=K-1$, and subsequently to others. Finally, substituting the estimation of $\mathbb{E}[\prod_{k=1}^K(\left\|\boldsymbol{X}_k\right\|_2+\varepsilon)]$ into \eqref{eq: law of total expectation} yields
	\begin{equation}\label{eq: Y prod estimate}
		\mathbb{E}\left[\left\|\boldsymbol{Y}_K  \cdots \boldsymbol{Y}_1\right\|_2\right] \leq \rho\gamma^{K},
	\end{equation}
	where $\gamma=\text{G}(\gamma_1,\ldots,\gamma_K)\left(1-\frac{2 \left(\underline{u} \lambda_d\left(\boldsymbol{A}^{\top}\boldsymbol{A}\right)-\overline{u} \sigma\right)}{ \overline{u}\left(\lambda_1\left(\boldsymbol{A}^{\top}\boldsymbol{A}\right)+\lambda_d\left(\boldsymbol{A}^{\top}\boldsymbol{A}\right)\right)}+\varepsilon\right)$ and $\text{G}(\gamma_1,\ldots,\gamma_K)=\left(\prod\limits_{k=1}^{K}\gamma_k\right)^{\frac{1}{K}}$. Combined with \eqref{eq: norm inequality}, we obtain \eqref{eq: linear convergence rate initial}-\eqref{eq: gamma}.
\end{proof}

It is worth noting that \eqref{eq: linear convergence rate initial}-\eqref{eq: gamma} implies a theoretical range of $\varepsilon$ to maintain an R-linear convergence rate of \hyperref[alg:SGA-RMSProp]{SGA-RMSProp} on the consistent LLSP, that is, 
\begin{equation*}
	0< \varepsilon < \frac{2 \left(\underline{u} \lambda_d\left(\boldsymbol{A}^{\top}\boldsymbol{A}\right)-\overline{u} \sigma\right)}{ \overline{u}\left(\lambda_1\left(\boldsymbol{A}^{\top}\boldsymbol{A}\right)+\lambda_d\left(\boldsymbol{A}^{\top}\boldsymbol{A}\right)\right)}.
\end{equation*}
Due to the impact of multiplier $\text{G}(\gamma_1,\ldots,\gamma_K)\leq1$, this range may become larger. Numerical experiments on the selection of $\varepsilon$ are presented in \cref{subsubsec: selection of epsilon}. Moreover, $\sigma$ can be viewed as an upper bound on the matrix analogue of the standard deviation of $\boldsymbol{M}_k$, $k=1,\ldots,K$. According to \eqref{eq: sigma}, it monotonically decreases to 0 when the batch size $B$ increases, making the upper bound of $\varepsilon$ positive for large $B$. As a result, it is straightforward to determine the values of $B$ and $\varepsilon$ that yield the R-linear convergence rate from \eqref{eq: linear convergence rate initial}-\eqref{eq: gamma}. This is shown in the corollary below, which is also a special case of \cref{thm: convergence of consistent LLSP}.
\begin{corollary}\label{thm: convergence of consistent LLSP corollary}
	Setting $B\geq\frac{4\overline{u}^2\log(2d)\max\limits_{\scriptscriptstyle j=1,\ldots,n}\left\{\|\boldsymbol{a}_j\|_{2}^{2}/p_j\right\}}{\underline{u}^2}\left(\frac{2\sqrt{2}\|\boldsymbol{A}\|_2}{\lambda_d\left(\boldsymbol{A}^{\top}\boldsymbol{A}\right)}+\left(\frac{\underline{u}}{3\overline{u}\lambda_d\left(\boldsymbol{A}^{\top}\boldsymbol{A}\right)}\right)^{\frac{1}{2}}\right)^2$ ensures that $\sigma\leq\frac{\underline{u} \lambda_d\left(\boldsymbol{A}^{\top}\boldsymbol{A}\right)}{4\overline{u}}$. Combined with setting $0<\varepsilon \leq \frac{\underline{u} \lambda_d\left(\boldsymbol{A}^{\top}\boldsymbol{A}\right)}{2\overline{u}\left(\lambda_1\left(\boldsymbol{A}^{\top}\boldsymbol{A}\right)+\lambda_d\left(\boldsymbol{A}^{\top}\boldsymbol{A}\right)\right)}$, we have
	\begin{equation}\label{eq: linear convergence rate}
		\textstyle\mathbb{E}[\left\|\boldsymbol{x}_{K+1}-\boldsymbol{x}^*\right\|_2] \leq  \rho\left\|\boldsymbol{x}_1-\boldsymbol{x}^*\right\|_2\left(G(\gamma_1,\ldots,\gamma_K)\left(1- \frac{\underline{u} \lambda_d\left(\boldsymbol{A}^{\top}\boldsymbol{A}\right)}{\overline{u}\left(\lambda_1\left(\boldsymbol{A}^{\top}\boldsymbol{A}\right)+\lambda_d\left(\boldsymbol{A}^{\top}\boldsymbol{A}\right)\right)} \right)\right)^K,
	\end{equation}
	in which $G(\gamma_1,\ldots,\gamma_K)\left(1- \frac{\underline{u} \lambda_d\left(\boldsymbol{A}^{\top}\boldsymbol{A}\right)}{\overline{u}\left(\lambda_1\left(\boldsymbol{A}^{\top}\boldsymbol{A}\right)+\lambda_d\left(\boldsymbol{A}^{\top}\boldsymbol{A}\right)\right)}\right)$ is a positive constant smaller than $1$.
\end{corollary}

\begin{proof}
	Firstly, by \eqref{eq: gamma}, setting $0<\sigma\leq\frac{\underline{u} \lambda_d\left(\boldsymbol{A}^{\top}\boldsymbol{A}\right)}{4\overline{u}}$ and $0<\varepsilon \leq \frac{\underline{u} \lambda_d\left(\boldsymbol{A}^{\top}\boldsymbol{A}\right)}{2\overline{u}\left(\lambda_1\left(\boldsymbol{A}^{\top}\boldsymbol{A}\right)+\lambda_d\left(\boldsymbol{A}^{\top}\boldsymbol{A}\right)\right)}$ gives \eqref{eq: linear convergence rate} directly. Next, let $c_1:=(2(\max\limits_{j=1,\ldots,n}\left\{\|\boldsymbol{a}_j\|_{2}^{2}/p_j\right\}-\lambda_d(\boldsymbol{A}^{\top}\boldsymbol{A}))\|\boldsymbol{A}\|_2^2 \log (2 d))^{\frac{1}{2}}$, $c_2:=\frac{\log (2 d)}{3}(\max\limits_{j=1,\ldots,n}\left\{\|\boldsymbol{a}_j\|_{2}^{2}/p_j\right\}-\lambda_d (\boldsymbol{A}^{\top} \boldsymbol{A})) $. Then \eqref{eq: sigma} implies that $-\sigma(\sqrt{B})^2+c_1 \sqrt{B}+c_2 = 0$, so $B =\left(\frac{c_1+\sqrt{c_1^2+4 \sigma c_2}}{2 \sigma}\right)^2$. In order to obtain $\sigma\leq\frac{\underline{u} \lambda_d\left(\boldsymbol{A}^{\top}\boldsymbol{A}\right)}{4\overline{u}}$, we need

	\begin{equation}\label{eq: B esitimate ori}\textstyle
		B \geq \bigg(2\overline{u}\bigg(c_1+\sqrt{c_1^2+\frac{\underline{u} c_2 \lambda_d\left(\boldsymbol{A}^{\top}\boldsymbol{A}\right)}{\overline{u}}}\bigg)\bigg/\left(\underline{u} \lambda_d\left(\boldsymbol{A}^{\top}\boldsymbol{A}\right)\right)\bigg)^2. 
	\end{equation}
	The following is the estimation of $c_1+\sqrt{c_1^2+\frac{\underline{u} c_2 \lambda_d\left(\boldsymbol{A}^{\top}\boldsymbol{A}\right)}{\overline{u}}}$:
	\begin{align*}\textstyle
		&c_1+\sqrt{c_1^2+\frac{\underline{u} c_2 \lambda_d\left(\boldsymbol{A}^{\top}\boldsymbol{A}\right)}{\overline{u}}}\nonumber \\
		\textstyle=&\textstyle \left(2\left(\max\limits_{j=1,\ldots,n}\left\{\frac{\|\boldsymbol{a}_j\|_{2}^{2}}{p_j}\right\}-\lambda_d\left(\boldsymbol{A}^{\top}\boldsymbol{A}\right)\right)\|\boldsymbol{A}\|_2^2 \log (2 d)\right)^{\frac{1}{2}}\nonumber \\
		\textstyle&\textstyle+\!\! \Bigg(\!2\!\left(\max\limits_{j=1,\ldots,n}\!\left\{\!\frac{\|\boldsymbol{a}_j\|_{2}^{2}}{p_j}\!\right\}\!-\!\lambda_d\!\left(\!\boldsymbol{A}^{\!\top}\!\boldsymbol{A}\!\right)\!\!\right)\!\|\boldsymbol{A}\|_2^2 \log (2 d)\! +\!\frac{\underline{u} \lambda_d\!\left(\!\boldsymbol{A}^{\!\top}\!\boldsymbol{A}\!\right)\!\log(2d)}{3 \overline{u}}\!\left(\max\limits_{j=1,\ldots,n}\!\left\{\!\frac{\|\boldsymbol{a}_j\|_{2}^{2}}{p_j}\!\right\}\!-\!\lambda_d\! \left(\!\boldsymbol{A}^{\!\top}\! \boldsymbol{A}\!\right)\!\!\right)\!\!\Bigg)^{\!\!\frac{1}{2}} \nonumber \\
		\textstyle\leq&\textstyle\left(\!2\!\!\max\limits_{j=1,\ldots,n}\!\!\left\{\!\frac{\|\boldsymbol{a}_j\|_{2}^{2}}{p_j}\!\right\}\!\|\boldsymbol{A}\|_2^2 \log (2 d)\!\!\right)^{\!\!\frac{1}{2}}\!\! +\!\Bigg(\!2\!\!\max\limits_{j=1,\ldots,n}\!\!\left\{\!\frac{\|\boldsymbol{a}_j\|_{2}^{2}}{p_j}\!\right\}\!\|\boldsymbol{A}\|_2^2\! \log (2 d)\!+\! \frac{\max\limits_{j=1,\ldots,n}\!\!\left\{\!\frac{\|\boldsymbol{a}_j\|_{2}^{2}}{p_j}\!\right\}\!\log(2d)\underline{u} \lambda_d\!\left(\!\boldsymbol{A}^{\!\top}\!\boldsymbol{A}\!\right)}{3 \overline{u}}\!\Bigg)^{\!\!\frac{1}{2}} \nonumber \\
		\textstyle=&\textstyle\left(\max\limits_{j=1,\ldots,n}\left\{\frac{\|\boldsymbol{a}_j\|_{2}^{2}}{p_j}\right\}\log(2d)\right)^{\frac{1}{2}}\left(\sqrt{2}\|\boldsymbol{A}\|_2+\left(2\|\boldsymbol{A}\|_2^2+\frac{\underline{u} \lambda_d(\boldsymbol{A}^{\top}\boldsymbol{A})}{3\overline{u}}\right)^{\frac{1}{2}}\right)\nonumber \\
		\textstyle\leq&\textstyle\left(\max\limits_{j=1,\ldots,n}\left\{\frac{\|\boldsymbol{a}_j\|_{2}^{2}}{p_j}\right\}\log(2d)\right)^{\frac{1}{2}}\left(2\sqrt{2}\|\boldsymbol{A}\|_2+\left(\frac{\underline{u} \lambda_d\left(\boldsymbol{A}^{\top}\boldsymbol{A}\right)}{3\overline{u}}\right)^{\frac{1}{2}}\right),
	\end{align*}
	where the second inequality is due to $\left(y^2+z^2\right)^{\frac{1}{2}}\leq |y|+|z|$. Combining the above with \eqref{eq: B esitimate ori}, we conclude that $\sigma\leq\frac{\underline{u} \lambda_d\left(\boldsymbol{A}^{\top}\boldsymbol{A}\right)}{4\overline{u}}$ follows from
	\begin{equation*}
		B\geq\frac{4\overline{u}^2\log(2d)\max\limits_{j=1,\ldots,n}\left\{\frac{\|\boldsymbol{a}_j\|_{2}^{2}}{p_j}\right\}}{\underline{u}^2}\left(\frac{2\sqrt{2}\|\boldsymbol{A}\|_2}{\lambda_d\left(\boldsymbol{A}^{\top}\boldsymbol{A}\right)}+\left(\frac{\underline{u}}{3\overline{u}\lambda_d\left(\boldsymbol{A}^{\top}\boldsymbol{A}\right)}\right)^{\frac{1}{2}}\right)^2.
	\end{equation*}
\end{proof}

\begin{remark}\label{remark:probability}
	A straightforward calculation shows that the minimum of $\max\limits_{j=1,\ldots,n}\{\|\boldsymbol{a}_j\|_2^2/p_j\}$ over all sample probabilities $\{p_j\}$ with $p_j>0$, $j=1,\ldots,n$, is attained when sampling rows with probability proportional to their squared length, that is, $p_j=\|\boldsymbol{a}_j\|_2^2/\|\boldsymbol{A}\|_{\mathrm{F}}^2$, $j=1,\ldots,n$. Under this choice, the resulting minimum is $\|\boldsymbol{A}\|_{\mathrm{F}}^2$. A similar statement can be found in \cite{needell2016stochastic}.
\end{remark}

\subsection{Inconsistent LLSP}
For the inconsistent LLSP, $\boldsymbol{A}\boldsymbol{x}^*$ is not equal to $\boldsymbol{b}$, meaning that $\nabla f_i(\boldsymbol{x}^*)\neq \boldsymbol{0}$ for some $i = 1,\ldots,n$. In this case, an extra residual term will appear in the upper bound on $\|\boldsymbol{x}_{K+1}-\boldsymbol{x}^*\|_2$ given by \eqref{eq: norm inequality}. This requires additional estimation to analyze the convergence of \hyperref[alg:SGA-RMSProp]{SGA-RMSProp}, that is, to show that the distances between iteration points $\{\boldsymbol{x}_k\}$ and a neighborhood of $\boldsymbol{x}^*$ converge to zero. The following are details.

Denote $\boldsymbol{r}^*=(r_1^*,\ldots,r_d^*)^{\top}:=\boldsymbol{A} \boldsymbol{x}^*-\boldsymbol{b}$, that is, $r_i^*:=\boldsymbol{a}_i^{\top} \boldsymbol{x}^*-b_i$, $i=1,\ldots,n$. Then the stochastic gradient at $\boldsymbol{x}_k$ can be written as
\begin{align}\label{eq:g_k inconsistent}
	\boldsymbol{g}_k &=\frac{1}{B} \sum_{i=1}^{B} \left(\frac{1}{p_{\scriptscriptstyle \xi_i^{(k)}}\!} \boldsymbol{a}_{\scriptscriptstyle \xi_i^{(k)}}\left(\boldsymbol{a}_{\scriptscriptstyle \xi_i^{(k)}}^{\top} \boldsymbol{x}_k-b_{\scriptscriptstyle \xi_i^{(k)}}\right)\right)\nonumber\\
	& =\frac{1}{B} \sum_{i=1}^{B}\left(\frac{1}{p_{\scriptscriptstyle \xi_i^{(k)}}\!} \boldsymbol{a}_{\scriptscriptstyle \xi_i^{(k)}}\left(\boldsymbol{a}_{\scriptscriptstyle \xi_i^{(k)}}^{\top} \boldsymbol{x}_k-\boldsymbol{a}_{\scriptscriptstyle \xi_i^{(k)}}^{\top} \boldsymbol{x}^*+r_{\scriptscriptstyle \xi_i^{(k)}}^*\right)\right) \nonumber\\
	& =\frac{1}{B} \sum_{i=1}^{B} \left(\frac{1}{p_{\scriptscriptstyle \xi_i^{(k)}}\!} \boldsymbol{a}_{\scriptscriptstyle \xi_i^{(k)}} \boldsymbol{a}_{\scriptscriptstyle \xi_i^{(k)}}^{\top}\left(\boldsymbol{x}_k-\boldsymbol{x}^*\right)+\frac{r_{\scriptscriptstyle \xi_i^{(k)}}^*}{p_{\scriptscriptstyle \xi_i^{(k)}}} \boldsymbol{a}_{\scriptscriptstyle \xi_i^{(k)}}\right) \nonumber \\
	& =\boldsymbol{M}_k\left(\boldsymbol{x}_k-\boldsymbol{x}^*\right)+\boldsymbol{h}_k, 
\end{align}
where $\boldsymbol{h}_k:=\frac{1}{B} \sum_{i=1}^{B} \frac{r_{\scriptscriptstyle \xi_i^{(k)}}^*}{p_{\scriptscriptstyle \xi_i^{(k)}}} \boldsymbol{a}_{\scriptscriptstyle \xi_i^{(k)}}$. Below is a lemma that provides an estimation of \ $\mathbb{E}\left[\left\|\boldsymbol{h}_k\right\|_2\right]$, which is derived from the proof of Theorem 5 in \cite{bollapragada2023fast}.
\begin{lemma}\label{lem:h_k estimate}
	For each $k=1,\ldots,K$, we have
	\begin{equation}\label{eq:h_k estimate}
		\mathbb{E}\left[\left\|\boldsymbol{h}_{k}\right\|_2\right] \leq \vcenter{\hbox{$\displaystyle \sqrt{\frac{2 \max\limits_{j=1,\ldots,n}\left\{\frac{\|\boldsymbol{a}_j\|_{2}^{2}}{p_j}\right\}\left\|\boldsymbol{r}^*\right\|_2^2 \log (d+1)}{B}} +\frac{\max\limits_{j=1,\ldots,n} \left\{\frac{\left\|\boldsymbol{a}_j\right\|_2|r_j^{*}|}{p_j}\right\} \log (d+1)}{3 B}$.}}
	\end{equation}
\end{lemma}

Next, we present the following theorem to analyze the convergence of \hyperref[alg:SGA-RMSProp]{SGA-RMSProp} on the inconsistent LLSP.

\begin{theorem} \label{thm: convergence of inconsistent LLSP}
	Let $\boldsymbol{x}^*$ be the optimal point of the inconsistent LLSP \eqref{eq:LLSP obj}, $\boldsymbol{r}^*=\boldsymbol{A} \boldsymbol{x}^*-\boldsymbol{b}$, $0<\varepsilon \leq \frac{\underline{u} \lambda_d\left(\boldsymbol{A}^{\top}\boldsymbol{A}\right)}{2\overline{u}\left(\lambda_1\left(\boldsymbol{A}^{\top}\boldsymbol{A}\right)+\lambda_d\left(\boldsymbol{A}^{\top}\boldsymbol{A}\right)\right)}$, $B\geq\frac{4\overline{u}^2\log(2d)\max\limits_{\scriptscriptstyle j=1,\ldots,n}\left\{\|\boldsymbol{a}_j\|_{2}^{2}/p_j\right\}}{\underline{u}^2}\left(\frac{2\sqrt{2}\|\boldsymbol{A}\|_2}{\lambda_d\left(\boldsymbol{A}^{\top}\boldsymbol{A}\right)}+\left(\frac{\underline{u}}{3\overline{u}\lambda_d\left(\boldsymbol{A}^{\top}\boldsymbol{A}\right)}\right)^{\frac{1}{2}}\right)^2$, $\eta_k\equiv\eta=\frac{2}{\overline{u}\left(\lambda_1\left(\boldsymbol{A}^{\top}\boldsymbol{A}\right)+\lambda_d\left(\boldsymbol{A}^{\top}\boldsymbol{A}\right)\right)}$, and $\{\boldsymbol{x}_k\}_{k=1}^{K}$ be the sequence generated by \hyperref[alg:SGA-RMSProp]{SGA-RMSProp}. Then we have
	\begin{equation*}
		\mathbb{E}[\left\|\boldsymbol{x}_{K+1}-\boldsymbol{x}^*\right\|_2] \leq \rho \left\|\boldsymbol{x}_1-\boldsymbol{x}^*\right\|_2\left(G(\gamma_1,\ldots,\gamma_K)\left(1- \frac{\underline{u} \lambda_d\left(\boldsymbol{A}^{\top}\boldsymbol{A}\right)}{\overline{u}\left(\lambda_1\left(\boldsymbol{A}^{\top}\boldsymbol{A}\right)+\lambda_d\left(\boldsymbol{A}^{\top}\boldsymbol{A}\right)\right)} \right)\right)^K+R,
	\end{equation*}
	with $\rho$, $G(\gamma_1,\ldots,\gamma_K)$ defined as in \cref{thm: convergence of consistent LLSP} and
	\begin{equation}\label{eq: confusion region size}
		\textstyle
		R \leq \frac{2\overline{u}^{\frac{3}{2}}}{\underline{u}^{\frac{3}{2}}\lambda_d\left(\boldsymbol{A}^{\top}\boldsymbol{A}\right)}\Bigg(\vcenter{\hbox{$\sqrt{\frac{2 \max\limits_{j=1,\ldots,n}\left\{\frac{\|\boldsymbol{a}_j\|_{2}^{2}}{p_j}\right\}\|\boldsymbol{r}^*\|_2^2 \log (d+1)}{B}}+\frac{\max\limits_{j=1,\ldots,n} \left\{\frac{\left\|\boldsymbol{a}_j\right\|_2|r_j^{*}|}{p_j}\right\} \log (d+1)}{3 B}$}}\Bigg).
	\end{equation}
\end{theorem}
\begin{proof}
	Substituting \eqref{eq:g_k inconsistent} into \eqref{eq: vector form of x in RMSProp}, for each $k=1,\ldots,K$, we have
	\begin{align}
		\boldsymbol{x}_{k+1}-\boldsymbol{x}^* & =\boldsymbol{x}_k-\boldsymbol{x}^*-\eta \boldsymbol{D}_k \boldsymbol{M}_k\left(\boldsymbol{x}_k-\boldsymbol{x}^*\right)-\eta \boldsymbol{D}_k \boldsymbol{h}_k \nonumber \\
		& =\left(\boldsymbol{I}-\eta \boldsymbol{D}_k \boldsymbol{M}_k\right)\left(\boldsymbol{x}_k-\boldsymbol{x}^*\right)-\eta \boldsymbol{D}_k \boldsymbol{h}_k \nonumber \\
		& =\boldsymbol{Y}_k\left(\boldsymbol{x}_k-\boldsymbol{x}^*\right)-\eta \boldsymbol{D}_k \boldsymbol{h}_k
	\end{align}
	where, consistent with \cref{thm: convergence of consistent LLSP}, we denote $\boldsymbol{Y}_k=\boldsymbol{I}-\eta \boldsymbol{D}_k \boldsymbol{M}_k$. It follows that
	\begin{equation}\textstyle
		\boldsymbol{x}_{K+1}-\boldsymbol{x}^*=\left(\boldsymbol{Y}_K \cdots \boldsymbol{Y}_1\right)\left(\boldsymbol{x}_1-\boldsymbol{x}^*\right)-\eta \sum\limits_{j=1}^K\left(\left(\prod\limits_{i=j+1}^K \boldsymbol{Y}_i\right) \boldsymbol{D}_j \boldsymbol{h}_j\right),
	\end{equation}
	in which $\prod_{i=K+1}^K \boldsymbol{Y}_i$ represents the unit matrix $\boldsymbol{I}$. We use this notation just for simplicity.
	
	By triangle inequality, we have
	\begin{equation}\label{eq: thm2 X_k estimate}\textstyle
		\| \boldsymbol{x}_{K+1}-\boldsymbol{x}^*\|_2\leq \left\| \boldsymbol{Y}_K  \cdots \boldsymbol{Y}_1\right\|_2 \left\|\boldsymbol{x}_1-\boldsymbol{x}^*\right\|_2 + \eta \left\| \sum\limits_{j=1}^K\left(\left(\prod\limits_{i=j+1}^K \boldsymbol{Y}_i\right) \boldsymbol{D}_j \boldsymbol{h}_j\right) \right\|_2.
	\end{equation}
	The term $\| \boldsymbol{Y}_K \cdots \boldsymbol{Y}_1\|_2$ has been studied in \cref{thm: convergence of consistent LLSP}, so we focus on the second part. By \cref{lem:SGA-RMSProp control}, $\|\boldsymbol{D}_k^{\frac{1}{2}}\boldsymbol{D}_{k-1}^{-\frac{1}{2}}\|_2$ is bounded above by $1$ for $k=1,\ldots,K$, while $\|\boldsymbol{D}_K^{\frac{1}{2}}\|_2 \|\boldsymbol{D}_j^{-\frac{1}{2}}\|_2$ is bounded above by $\overline{u}^\frac{1}{2}/\underline{u}^{\frac{1}{2}}$ for $j=1,\ldots,K-1$. Similarly to \eqref{eq: Y transfer} and \eqref{eq: D_k to X_k}, we have
	\begin{equation}\label{eq: thm2 total Y inequality}
		\left\|\boldsymbol{Y}_K \boldsymbol{Y}_{K-1} \cdots \boldsymbol{Y}_{j+1}\right\|_2  \leq \frac{\overline{u}^\frac{1}{2}}{\underline{u}^{\frac{1}{2}}} \prod_{i=j+1}^K\left(\left\|\boldsymbol{X}_i\right\|_2+\varepsilon\right),~j=1,\ldots,K-1,
	\end{equation}
	and
	\begin{align*}
		\textstyle\left\|\sum\limits_{j=1}^K\!\!\left(\!\!\left(\!\prod\limits_{i=j+1}^K\!\! \boldsymbol{Y}_i\!\right)\!\! \boldsymbol{D}_j \boldsymbol{h}_j\!\!\right)\!\right\|_{\!2} &\textstyle\!\!\!\leq\!\sum\limits_{j=1}^K\!\left\|\!\left(\!\prod\limits_{i=j+1}^K\!\! \boldsymbol{Y}_i\!\right)\! \!\boldsymbol{D}_j \boldsymbol{h}_j\!\right\|_{\!2}  \\
		\textstyle&\textstyle\!\!\!\leq\! \sum\limits_{j=1}^K\!\!\left(\left\|\prod\limits_{i=j+1}^K\!\! \boldsymbol{Y}_i\right\|_{\!2}\!\!\!\left\|\boldsymbol{D}_j\right\|_{\!2}\! \left\| \boldsymbol{h}_j\right\|_{\!2}\!\!\right)  
		\!\!\overset{\eqref{eq: D_k estimate},\eqref{eq: thm2 total Y inequality}}{\leq}\!\! \frac{\overline{u}^{\frac{3}{2}}}{\underline{u}^{\frac{1}{2}}}\!  \sum\limits_{j=1}^K\!\!\left(\!\!\left(\!\prod\limits_{i=j+1}^K\!\!\left(\left\|\boldsymbol{X}_i\right\|_{\!2}\!+\!\varepsilon\right)\!\right)\!\left\| \boldsymbol{h}_j\right\|_{\!2}\!\!\right)\!,
	\end{align*}
	where $\boldsymbol{X}_i=\boldsymbol{I}-\eta \boldsymbol{D}_{i-1}^{\frac{1}{2}} \boldsymbol{M}_i \boldsymbol{D}_{i-1}^{\frac{1}{2}}$, $i =2,\ldots,K$, and the notation $\prod_{i=K+1}^K\left(\left\|\boldsymbol{X}_i\right\|_2+\varepsilon\right)$ represents $1$ for simplicity. Take expectation on both sides and we have
	\begin{equation}\label{eq: thm2 second part}\textstyle
		\mathbb{E}\left[\left\|\sum\limits_{j=1}^K\left(\left(\prod\limits_{i=j+1}^K \boldsymbol{Y}_i\right) \boldsymbol{D}_j \boldsymbol{h}_j\right)\right\|_2 \right] \leq \frac{\overline{u}^{\frac{3}{2}}}{\underline{u}^{\frac{1}{2}}}\sum\limits_{j=1}^K\left( \mathbb{E}\left[ \left(\prod\limits_{i=j+1}^K\left(\left\|\boldsymbol{X}_i\right\|_2+\varepsilon\right)\right)\left\| \boldsymbol{h}_j\right\|_2\right]\right).
	\end{equation}
	
	Recalling \cref{thm: convergence of consistent LLSP corollary}, setting $B$ and $\varepsilon$ as in \cref{thm: convergence of inconsistent LLSP} yields
	\begin{equation}\label{eq: linear convergence rate inequality}
		1- \frac{2 \left(\underline{u} \lambda_d\left(\boldsymbol{A}^{\top}\boldsymbol{A}\right)-\overline{u} \sigma\right)}{ \overline{u}\left(\lambda_1\left(\boldsymbol{A}^{\top}\boldsymbol{A}\right)+\lambda_d\left(\boldsymbol{A}^{\top}\boldsymbol{A}\right)\right)}+\varepsilon\leq 1- \frac{\underline{u} \lambda_d\left(\boldsymbol{A}^{\top}\boldsymbol{A}\right)}{\overline{u}\left(\lambda_1\left(\boldsymbol{A}^{\top}\boldsymbol{A}\right)+\lambda_d\left(\boldsymbol{A}^{\top}\boldsymbol{A}\right)\right)}.
	\end{equation}
	Notice that $\boldsymbol{h}_j$ is $\mathcal{F}_{j+1}$-measurable, $j=1,\ldots,K-1$. Applying the process similar to \eqref{eq: law of total expectation}, \eqref{eq: measurable expectation formula}, and \eqref{eq: X_k estimate}, $\mathbb{E}\left[ \left(\prod_{i=j+1}^K\left(\left\|\boldsymbol{X}_i\right\|_2+\varepsilon\right)\right)\left\| \boldsymbol{h}_j\right\|_2\right]$ can be estimated as follows:
	\begin{align}\label{eq: thm2 prod estimate}
		&\textstyle\mathbb{E}\left[ \left(\prod\limits_{i=j+1}^K\left(\left\|\boldsymbol{X}_i\right\|_2+\varepsilon\right)\right)\left\| \boldsymbol{h}_j\right\|_2\right] \nonumber \\
		\textstyle\overset{\eqref{eq: law of total expectation},\eqref{eq: measurable expectation formula}}{=}\hspace{1.6mm}&\textstyle\mathbb{E}\left[\cdots\mathbb{E}\left[\mathbb{E}\left[\left\|\boldsymbol{X}_K\right\|_2+\varepsilon|\mathcal{F}_K\right]\left(\prod\limits_{i=j+1}^{K-1}\left(\left\|\boldsymbol{X}_i\right\|_2+\varepsilon\right)\right)\left\| \boldsymbol{h}_j\right\|_2\Bigg|\mathcal{F}_{K-1}\right]\cdots \right]\nonumber \\
		\overset{\eqref{eq: X_k estimate}, \eqref{eq: linear convergence rate inequality}}{\leq}& \left(1- \frac{\underline{u} \lambda_d\left(\boldsymbol{A}^{\top}\boldsymbol{A}\right)}{\overline{u}\left(\lambda_1\left(\boldsymbol{A}^{\top}\boldsymbol{A}\right)+\lambda_d\left(\boldsymbol{A}^{\top}\boldsymbol{A}\right)\right)} \right) \nonumber \\
		\textstyle&\textstyle\mathbb{E}\left[\cdots \mathbb{E}\left[\mathbb{E}\left[\left\|\boldsymbol{X}_{K-1}\right\|_2+\varepsilon|\mathcal{F}_{K-1}\right]\left(\prod\limits_{i=j+1}^{K-2}\left(\left\|\boldsymbol{X}_i\right\|_2+\varepsilon\right)\right)\left\| \boldsymbol{h}_j\right\|_2\Bigg|\mathcal{F}_{K-2}\right]\cdots \right] \nonumber \\
		\textstyle\leq\hspace{6.2mm}&\textstyle \ldots \nonumber \\
		\textstyle\leq\hspace{6.2mm}&\textstyle \left(1- \frac{\underline{u} \lambda_d\left(\boldsymbol{A}^{\top}\boldsymbol{A}\right)}{\overline{u}\left(\lambda_1\left(\boldsymbol{A}^{\top}\boldsymbol{A}\right)+\lambda_d\left(\boldsymbol{A}^{\top}\boldsymbol{A}\right)\right)} \right)^{K-j-1}\mathbb{E}\left[\mathbb{E}\left[\left\|\boldsymbol{X}_{j+1}\right\|_2+\varepsilon|\mathcal{F}_{j+1}\right]\left\| \boldsymbol{h}_j\right\|_2\right] \nonumber \\
		\textstyle\leq\hspace{6.2mm}&\textstyle \left(1- \frac{\underline{u} \lambda_d\left(\boldsymbol{A}^{\top}\boldsymbol{A}\right)}{\overline{u}\left(\lambda_1\left(\boldsymbol{A}^{\top}\boldsymbol{A}\right)+\lambda_d\left(\boldsymbol{A}^{\top}\boldsymbol{A}\right)\right)}\right)^{K-j} \mathbb{E}\left[\left\| \boldsymbol{h}_j\right\|_2\right]
	\end{align}
	for $j=1,\ldots, K-1$. Moreover, the above also holds for $j=K$ since $\prod_{i=K+1}^K\left(\left\|\boldsymbol{X}_i\right\|_2+\varepsilon\right)$ represents $1$.
	
	As mentioned in \cref{lem:h_k estimate}, $\mathbb{E}[\left\|\boldsymbol{h}_j\right\|_2]$ has the same value for each $j=1,\ldots,K$. Combining \eqref{eq: thm2 second part} with \eqref{eq: thm2 prod estimate}, we obtain
	\begin{equation}\label{eq: thm2 second part estimate}
		\textstyle
		\mathbb{E}\left[\left\|\sum \limits_{j=1}^K\left(\left(\prod \limits_{i=j+1}^K \boldsymbol{Y}_i\right) \boldsymbol{D}_j \boldsymbol{h}_j\right)\right\|_2 \right]
		\leq\frac{\overline{u}^{\frac{3}{2}} \mathbb{E}\left[\left\| \boldsymbol{h}_1\right\|_2\right]}{\underline{u}^{\frac{1}{2}}}\sum \limits_{j=1}^K \left(1- \frac{\underline{u} \lambda_d\left(\boldsymbol{A}^{\top}\boldsymbol{A}\right)}{\overline{u}\left(\lambda_1\left(\boldsymbol{A}^{\top}\boldsymbol{A}\right)+\lambda_d\left(\boldsymbol{A}^{\top}\boldsymbol{A}\right)\right)}\right)^{K-j}.
	\end{equation}
	
	Finally, taking expectation on both sides of \eqref{eq: thm2 X_k estimate} gives
	\begin{align}
		\textstyle\mathbb{E}[\left\|\boldsymbol{x}_{K+1}-\boldsymbol{x}^*\right\|_2] \overset{\eqref{eq: Y prod estimate}}{\leq}&\textstyle \rho \left\|\boldsymbol{x}_1-\boldsymbol{x}^*\right\|_2\left(G(\gamma_1,\ldots,\gamma_K)\left(1- \frac{\underline{u} \lambda_d\left(\boldsymbol{A}^{\top}\boldsymbol{A}\right)}{\overline{u}\left(\lambda_1\left(\boldsymbol{A}^{\top}\boldsymbol{A}\right)+\lambda_d\left(\boldsymbol{A}^{\top}\boldsymbol{A}\right)\right)}\right) \right)^K \nonumber \\
		&\textstyle+\eta\mathbb{E}\left[ \left\| \sum\limits_{j=1}^K\left(\left(\prod\limits_{i=j+1}^K \boldsymbol{Y}_i\right) \boldsymbol{D}_j \boldsymbol{h}_j\right) \right\|_2\right].\nonumber
	\end{align}
	
	According to \cref{lem:h_k estimate} and our setting of $\eta$, we have
	\begin{align*}
		R &\textstyle\hspace{2mm}:= \hspace{2mm}\eta\mathbb{E}\left[ \left\| \sum\limits_{j=1}^K\left(\left(\prod\limits_{i=j+1}^K \boldsymbol{Y}_i\right) \boldsymbol{D}_j \boldsymbol{h}_j\right) \right\|_2\right] \nonumber \\
		&\textstyle\overset{\eqref{eq: thm2 second part estimate}}{\leq}\frac{\eta \overline{u}^{\frac{3}{2}} \mathbb{E}\left[\left\| \boldsymbol{h}_1\right\|_2\right]}{\underline{u}^{\frac{1}{2}}}\sum\limits_{j=1}^K \left(1- \frac{\underline{u} \lambda_d\left(\boldsymbol{A}^{\top}\boldsymbol{A}\right)}{\overline{u}\left(\lambda_1\left(\boldsymbol{A}^{\top}\boldsymbol{A}\right)+\lambda_d\left(\boldsymbol{A}^{\top}\boldsymbol{A}\right)\right)}\right)^{K-j} \nonumber \\
		&\hspace{2.3mm}\leq \hspace{2.3mm}\frac{\eta \overline{u}^{\frac{5}{2}}\left(\lambda_1\left(\boldsymbol{A}^{\top}\boldsymbol{A}\right)+\lambda_d\left(\boldsymbol{A}^{\top}\boldsymbol{A}\right)\right) \mathbb{E}\left[\left\| \boldsymbol{h}_1\right\|_2\right]}{\underline{u}^{\frac{3}{2}}\lambda_d\left(\boldsymbol{A}^{\top}\boldsymbol{A}\right)}\nonumber \\ 
		&\textstyle\overset{\eqref{eq:h_k estimate}}{\leq}\frac{2\overline{u}^{\frac{3}{2}}}{\underline{u}^{\frac{3}{2}}\lambda_d\left(\boldsymbol{A}^{\top}\boldsymbol{A}\right)}\Bigg( \vcenter{\hbox{$\sqrt{\frac{2 \max\limits_{j=1,\ldots,n}\left\{\frac{\|\boldsymbol{a}_j\|_{2}^{2}}{p_j}\right\}\|\boldsymbol{r}^*\|_2^2 \log (d+1)}{B}}+\frac{\max\limits_{j=1,\ldots,n} \left\{\frac{\left\|\boldsymbol{a}_j\right\|_2|r_j^{*}|}{p_j}\right\} \log (d+1)}{3 B}$}}\Bigg).
	\end{align*}
	
	This completes the proof.
\end{proof}

Extending from the consistent case, \cref{thm: convergence of inconsistent LLSP} indicates that for the inconsistent LLSP, the sequence $\{\boldsymbol{x}_k\}$ generated by \hyperref[alg:SGA-RMSProp]{SGA-RMSProp} still converges R-linearly, but to a neighborhood of $\boldsymbol{x}^*$ with radius $R$. This kind of neighborhood can be viewed as the ``region of confusion'', where the method fails to obtain a clear direction towards the optimal point \cite{bertsekas2015convex}. As shown in \eqref{eq: confusion region size}, the region of confusion provided by \cref{thm: convergence of inconsistent LLSP} can be controlled by the batch size $B$. An increase in the batch size results in a reduction of the region of confusion.

\section{Numerical experiments}\label{sec:numerical experiments}
\renewcommand{\arraystretch}{0.6}
\hspace{-1mm}In this section, \hspace{-0.2mm}we evaluate the \hspace{-0.2mm}performance of \hyperref[alg:SGA-RMSProp]{SGA-RMSProp} on LLSP using synthetic and real data. We compare \hyperref[alg:SGA-RMSProp]{SGA-RMSProp} with SGD, showing that \hyperref[alg:SGA-RMSProp]{SGA-RMSProp} generally converges faster on LLSP with the small batch size and exhibits a faster initial convergence rate with the large batch size. Based on this, a strategy for switching from \hyperref[alg:SGA-RMSProp]{SGA-RMSProp} to SGD is proposed, combining the benefits of these two algorithms. Our experiments are performed in MATLAB R2023b (version 23.2.0.2365128) on a desktop equipped with 64 GB memory, an Intel Core i9-12900K (3.2GHz).

When sampling the mini-batch stochastic gradient, $B=50$ and $1000$ are used to indicate the small and large batch sizes respectively, and let $p_i=\frac{\|\boldsymbol{a}_i\|_2^2}{\|\boldsymbol{A}\|_{\mathrm{F}}^2}$ (see \cref{remark:probability}), $i=1,\ldots,n$. For parameter settings of \hyperref[alg:SGA-RMSProp]{SGA-RMSProp}, let $\beta_k = \tfrac{1}{2}\max \limits_{j=1, \ldots, d}\left\{\tfrac{g_{k, j}^2-\tfrac{1}{\underline{u}^2}}{g_{k, j}^2-u_{k-1, j}}, \frac{g_{k, j}^2-(1+\varepsilon)^2 u_{k-1, j}}{g_{k, j}^2-u_{k-1, j}},0\right\}+\frac{1}{2}$ when $\frac{g_{k,j}^2}{u_{k-1,j}}> 1$, $ j=1,\ldots,d$, within the range defined in \hyperref[alg:SGA-RMSProp]{SGA-RMSProp}. We use the following step size when $B=50$:
\begin{equation*}
	\eta= \begin{cases}\frac{1.1B}{\overline{u}\left(\|\boldsymbol{A}\|_{\mathrm{F}}^2+(B-1) \lambda_d(\boldsymbol{A}^{\top}\boldsymbol{A})\right)}, & \|\boldsymbol{A}\|_{\mathrm{F}}^2-(B-1)\left(\lambda_1(\boldsymbol{A}^{\top}\boldsymbol{A})-\lambda_d(\boldsymbol{A}^{\top}\boldsymbol{A})\right) \geq 0, \\ \frac{\left(2.1+\sqrt{\lambda_d(\boldsymbol{A}^{\top}\boldsymbol{A})/\lambda_1(\boldsymbol{A}^{\top}\boldsymbol{A})}\right) B}{\overline{u}\left(\|\boldsymbol{A}\|_{\mathrm{F}}^2+(B-1)\left(\lambda_1(\boldsymbol{A}^{\top}\boldsymbol{A})+\lambda_d(\boldsymbol{A}^{\top}\boldsymbol{A})\right)\right)}, & \|\boldsymbol{A}\|_{\mathrm{F}}^2-(B-1)\left(\lambda_1(\boldsymbol{A}^{\top}\boldsymbol{A})-\lambda_d(\boldsymbol{A}^{\top}\boldsymbol{A})\right)<0,\end{cases}
\end{equation*}
which is inspired by the step size of mini-batch SGD suggested by Moorman \cite{moorman2021randomized}.
When $B=1000$, we apply a heuristic step size based on our theoretical results. Experiments show that the algorithm always converges with $\eta=\frac{2}{\overline{u}(\lambda_{1}(\boldsymbol{A}^{\top}\boldsymbol{A})+\lambda_{d}(\boldsymbol{A}^{\top}\boldsymbol{A}))}$, which is also the step size in \cref{thm: convergence of consistent LLSP} and \cref{thm: convergence of inconsistent LLSP}, and a slightly larger $\eta$ may speed up convergence. Therefore, we set  $\eta=\frac{2+\sqrt{\lambda_d(\boldsymbol{A}^{\top}\boldsymbol{A})/\lambda_1(\boldsymbol{A}^{\top}\boldsymbol{A})}}{\overline{u}(\lambda_{1}(\boldsymbol{A}^{\top}\boldsymbol{A})+\lambda_{d}(\boldsymbol{A}^{\top}\boldsymbol{A}))}$ when $B=1000$. The effect of $\overline{u}$, $\underline{u}$, and $\varepsilon$ is discussed in \cref{subsubsec: selection of epsilon}. 

In each experiment, the algorithm will be run 100 times. Unless otherwise specified, we regard the algorithm as converged if $\frac{\|\boldsymbol{x}_k-\boldsymbol{x}^*\|_2}{\|\boldsymbol{x}_1-\boldsymbol{x}^*\|_2}\leq 10^{-4}$. When presenting the results, we calculate the mean using values from the 5th to 95th percentile to avoid extreme cases skewing the average.

\subsection{Experiments for consistent LLSP}
\label{subsec:experiments for consistent LLSP}
In this subsection, we present the convergence results of \hyperref[alg:SGA-RMSProp]{SGA-RMSProp} for different batch sizes, showing the effect of $\overline{u}$, $\underline{u}$, and $\varepsilon$. Next, the performances of \hyperref[alg:SGA-RMSProp]{SGA-RMSProp} and other algorithms are compared on the consistent LLSP.

\subsubsection{Data generation}\label{subsubsec:Data description}
We generate the synthetic data similar to \cite{strakovs1991real,bollapragada2023fast}. Let $\boldsymbol{A}\in\mathbb{R}^{n\times d}$ with $n=10^6$ and $d=10^2$. Its thin singular value decomposition is $\boldsymbol{U} \boldsymbol{\Sigma} \boldsymbol{V}^{\top}$, where $\boldsymbol{U} \in \mathbb{R}^{n \times d}$ has orthonormal columns, $\boldsymbol{V}$ is orthogonal, both chosen uniformly at random as in Proposition 7.2 of \cite{eaton1983multivariate}, and $\boldsymbol{\Sigma}=\text{diag}(s_1,\ldots,s_d)$. Then the singular values of $\boldsymbol{A}^{\top}\boldsymbol{A}$ are $\{s_1^2,\ldots,s_d^2\}$, set to exponential or algebraic decay. Specifically, given the condition number $\kappa=\frac{\lambda_1\left(\boldsymbol{A}^{\top}\boldsymbol{A}\right)}{\lambda_d\left(\boldsymbol{A}^{\top}\boldsymbol{A}\right)}$ of $\boldsymbol{A}^{\top}\boldsymbol{A}$, decay rate $q$, and $\lambda_d(\boldsymbol{A}^{\top}\boldsymbol{A})$, the formula for the exponential decay (ED) is:
\begin{equation*}
	s_j=\sqrt{\lambda_d(\boldsymbol{A}^{\top}\boldsymbol{A})+\left(\frac{d-j}{d-1}\right)\lambda_d(\boldsymbol{A}^{\top}\boldsymbol{A})(\kappa-1) q^{j-1}},~j=1, \ldots, d,
\end{equation*}
and the formula for the algebraic decay (AD) is
\begin{equation*}
	s_j=\sqrt{\lambda_d(\boldsymbol{A}^{\top}\boldsymbol{A})+\left(\frac{d-j}{d-1}\right)^q\lambda_d(\boldsymbol{A}^{\top}\boldsymbol{A})(\kappa-1)},~j=1,\ldots, d.
\end{equation*}
With the same condition number, different decay types and rates result in different singular value distributions \cite{strakovs1991real}. We take $\lambda_d(\boldsymbol{A}^{\top}\boldsymbol{A})=1$, $\kappa=20, 50, 100$ (small, medium, large condition numbers) and $q=0.2, 0.7$ for ED, $q=1, 2$ for AD. These problems are denoted as (decay-type, $\kappa$, $q$). Lastly, let the minimizer $\boldsymbol{x}^*$ have independent standard normal entries and $\boldsymbol{b}$ is set to $\boldsymbol{A}\boldsymbol{x}^*$. We have defined 12 problems, each with 3 random instances, totaling 36 instances. Since  $\boldsymbol{A}$ and $\boldsymbol{x}^*$ differ across the instances, the data is sufficiently randomized. Hence, we fix the initial point $\boldsymbol{x}_1 = (2,\ldots,2)^\top$ in all experiments. For each problem, we report the mean and standard deviation over three instances, in the form of (mean, standard deviation).
\begin{table}[H]\label{tab:results_different_u_epsilon1}
	\small
	\centering
	\caption{Numbers of iterations required for SGA-RMSProp to converge ($\varepsilon=\varepsilon_1$), reported as (mean, standard deviation) over three instances for $\overline{u}_i$, $\underline{u}_i$, $i=1,\ldots,4$.}
	\begin{tabular}{cccccc}
		\toprule
		Batch sizes &Problems & ($\underline{u}_1$, $\overline{u}_1$)& ($\underline{u}_2$, $\overline{u}_2$) & ($\underline{u}_3$, $\overline{u}_3$) &($\underline{u}_4$, $\overline{u}_4$) \\
		\midrule
		\multirow{12}{*}{$B=50$}&(ED, 20, 0.7) & (113, 0.6) & (115, 0.6) &(117, 0.0)&(121, 1.2) \\
		&(ED, 20, 0.2) & (113, 1.0) & (113, 0.6) & (116, 1.2)&(121, 2.0) \\
		&(ED, 50, 0.7) & (255, 2.5) & (257, 3.2) & (266, 3.0) &(275, 3.6) \\
		&(ED, 50, 0.2) & (261, 5.1) & (266, 2.0) & (276, 3.0)&(287, 3.0)  \\
		&(ED, 100, 0.7) & (493, 4.5) & (492, 2.0) & (506, 2.5)&(525, 2.3) \\
		&(ED, 100, 0.2) & (531, 13.3) & (550, 13.5) & (575, 18.3) &(600, 14.0) \\
		&(AD, 20, 1) & (149, 8.5) & (153, 9.5) & (157, 10.0) &(162, 11.7) \\
		&(AD, 20, 2) & (145, 2.1) & (142, 0.6) & (144, 0.6) & (147, 2.0)  \\
		&(AD, 50, 1) & (327, 41.0)  & (333, 39.9) & (341, 46.5) &(350, 49.4) \\
		&(AD, 50, 2) & (304, 5.7) & (311, 6.4) & (319, 7.2) &(327, 9.0) \\
		&(AD, 100, 1) & (541, 46.1) & (559, 53.4) & (569, 48.3) &(583, 59.3) \\
		&(AD, 100, 2) & (558, 39.2) & (565, 45.5) & (582, 43.8) &(595, 49.2)\\
		\midrule
		\multirow{12}{*}{$B=1000$}&(ED, 20, 0.7) & (97, 1.7) & (97, 3.1) & (101, 3.1) &(109, 4.0)\\
		&(ED, 20, 0.2) & (98, 2.1) & (100, 2.1) & (106, 2.5) & (119, 5.7) \\
		&(ED, 50, 0.7) & (233, 3.0) & (240, 3.1) & (251, 5.9) & (268, 4.9) \\
		&(ED, 50, 0.2) & (236, 1.5) & (244, 5.3) & (259, 9.5) & (281, 19.0) \\
		&(ED, 100, 0.7) & (459, 2.5) & (463, 3.5) & (477, 9.9) & (500, 13.6)\\
		&(ED, 100, 0.2) & (467, 5.6) & (479, 13.1) & (496, 18.9) & (530, 27.5)\\
		&(AD, 20, 1) & (209, 8.1) & (127, 5.1) & (80, 1.0) & (74, 3.5)\\
		&(AD, 20, 2) & (173, 20.0) & (110, 8.7) & (84, 5.3) & (78, 2.1) \\
		&(AD, 50, 1) & (254, 29.0) & (191, 18.4) & (172, 17.3) & (171, 21.0) \\
		&(AD, 50, 2) & (266, 27.9) & (206, 9.9) & (195, 4.0) & (195, 3.5)\\
		&(AD, 100, 1) & (366, 4.7) & (316, 16.9) & (298, 25.5) & (294, 33.8) \\
		&(AD, 100, 2) & (371, 21.9) & (360, 24.7) & (361, 23.1) &(362, 27.8) \\
		\bottomrule
	\end{tabular}
\end{table}

\subsubsection{Effect of $\overline{u}$, $\underline{u}$, and $\varepsilon$}
\label{subsubsec: selection of epsilon}
Values of $\overline{u}$ are chosen as $\overline{u}_i=(10^{-i}\min\limits_{j=1,\ldots,d}\{g_{1,j}^2|g_{1,j}\neq 0\})^{-\tfrac{1}{2}}$, $i=1,\ldots,4$. Note that $\boldsymbol{g}_1$ is directly computed in the initial step, with no parameters required. This setting allows the upper bound $\overline{u}$ and the initial value $\boldsymbol{u}_0$ to match the scale of the observed stochastic gradient. For each $\overline{u}_i$, we set $\underline{u}$ to $\underline{u}_i = \tfrac{\overline{u}_i}{5}$ and $\underline{u}_i^{\prime}=\tfrac{\overline{u}_i}{10}$, $i=1,\ldots,4$. Let $\varepsilon$ take three values, $\varepsilon_1=\frac{\lambda_{d}(\boldsymbol{A}^{\top}\boldsymbol{A})}{\lambda_{1}(\boldsymbol{A}^{\top}\boldsymbol{A})}$, $\varepsilon_2=\frac{5\lambda_{d}(\boldsymbol{A}^{\top}\boldsymbol{A})}{\lambda_{1}(\boldsymbol{A}^{\top}\boldsymbol{A})}$, $\varepsilon_3=\frac{10\lambda_{d}(\boldsymbol{A}^{\top}\boldsymbol{A})}{\lambda_{1}(\boldsymbol{A}^{\top}\boldsymbol{A})}$, representing low, medium, and high levels, respectively. The effect of these parameter combinations is tested in the following experiments.

\cref{tab:results_different_u_epsilon1} and \cref{tab:results_different_u_prime_epsilon1} show the numbers of iterations required for \hyperref[alg:SGA-RMSProp]{SGA-RMSProp} to converge with $\varepsilon=\varepsilon_1$, given different combinations of $\overline{u}$ and $\underline{u}$. The mean is rounded to the nearest integer and the standard deviation is rounded to the nearest tenth when presenting the results.

\vspace{-1.5mm}\begin{table}[h]\label{tab:results_different_u_prime_epsilon1}
	\small
	\centering
	\caption{Numbers of iterations required for SGA-RMSProp to converge ($\varepsilon=\varepsilon_1$), reported as (mean, standard deviation) over three instances for $\overline{u}_i$, $\underline{u}_i^{\prime}$, $i=1,\ldots,4$.}
	\vspace{-1.5mm}\begin{tabular}{cccccc}
		\toprule
		Batch sizes &Problems & ($\underline{u}_1^{\prime}$, $\overline{u}_1$)& ($\underline{u}_2^{\prime}$, $\overline{u}_2$) & ($\underline{u}_3^{\prime}$, $\overline{u}_3$) &($\underline{u}_4^{\prime}$, $\overline{u}_4$) \\
		\midrule
		\multirow{12}{*}{$B=50$}&(ED, 20, 0.7) & (114, 0.6) & (115, 0.0) &(117, 0.0)&(121, 1.2) \\
		&(ED, 20, 0.2) & (114, 0.6) & (113, 0.6) & (116, 1.2)&(121, 2.0) \\
		&(ED, 50, 0.7) & (256, 0.6) & (257, 1.5) & (266, 3.0) &(275, 3.6) \\
		&(ED, 50, 0.2) & (260, 3.8) & (262, 2.3) & (276, 3.0)&(287, 3.0)  \\
		&(ED, 100, 0.7) & (494, 4.4) & (489, 1.0) & (506, 2.5)&(525, 2.3) \\
		&(ED, 100, 0.2) & (528, 12.1) & (541, 9.8) & (575, 18.3) &(600, 14.0) \\
		&(AD, 20, 1) & (150, 9.3) & (152, 9.5) & (157, 10.0) &(162, 11.7) \\
		&(AD, 20, 2) & (146, 1.5) & (143, 1.5) & (144, 0.6) & (147, 2.0)  \\
		&(AD, 50, 1) & (329, 36.4)  & (327, 47.2) & (341, 46.5) &(350, 49.4) \\
		&(AD, 50, 2) & (305, 3.8) & (311, 6.2) & (319, 7.2) &(327, 9.0) \\
		&(AD, 100, 1) & (555, 33.0) & (559, 52.1) & (569, 48.3) &(583, 59.3) \\
		&(AD, 100, 2) & (558, 40.0) & (565, 43.2) & (582, 43.8) &(595, 49.2)\\
		\midrule
		\multirow{12}{*}{$B=1000$}&(ED, 20, 0.7) & (97, 1.7) & (97, 2.9) & (100, 3.2) &(110, 3.8)\\
		&(ED, 20, 0.2) & (98, 2.1) & (100, 2.0)  & (106, 2.5) & (118, 4.7) \\
		&(ED, 50, 0.7) & (233, 3.0) & (239, 7.0) & (251, 5.9) & (267, 4.6) \\
		&(ED, 50, 0.2) & (236, 1.5) & (245, 4.0) & (259, 9.1) & (281, 19.0) \\
		&(ED, 100, 0.7) & (459, 2.5) & (464, 5.3) & (477, 9.9) & (500, 13.6)\\
		&(ED, 100, 0.2) & (467, 5.6) & (481, 12.9) & (496, 18.9) & (530, 27.5)\\
		&(AD, 20, 1) & (209, 8.1) & (136, 15.9) & (80, 1.0) & (74, 3.5)\\
		&(AD, 20, 2) & (173, 20.0) & (112, 11.9) & (84, 5.3) & (78, 2.1) \\
		&(AD, 50, 1) & (254, 29.0) & (195, 16.2) & (172, 17.3) & (171, 21.0) \\
		&(AD, 50, 2) & (266, 27.9) & (211, 6.0) & (195, 4.0) & (195, 3.5)\\
		&(AD, 100, 1) & (366, 4.7) & (313, 12.5) & (298, 25.5) & (294, 33.8) \\
		&(AD, 100, 2) & (371, 21.9) & (360, 25.1) & (361, 23.1) &(362, 27.8) \\
		\bottomrule
	\end{tabular}
\end{table}

\vspace{-1.5mm}The results show that, under our parameters selection strategy where $\overline{u}$ is specified first and $\underline{u}$ is chosen accordingly, the algorithm exhibits similar performance when using the parameters $(\underline{u}_i, \overline{u}_i)$ and $(\underline{u}_i^{\prime}, \overline{u}_i)$ respectively, $i=1,\ldots,4$. Therefore, we use $\underline{u}_i = \tfrac{\overline{u}_i}{5}$, $i=1,\ldots,4$, in the following experiments.

\cref{tab:results_different_u_epsilon1}-\cref{tab:results_different_u_epsilon3} indicate that the algorithm generally performs better with $\overline{u}=\overline{u}_1$, the smallest value among the considered parameters, when $B=50$. According to \eqref{eq: D_k estimate} in \cref{lem:SGA-RMSProp control}, we have $\|\boldsymbol{D}_k\|_2\leq \overline{u}$, $k=1,\ldots,K$. Since the uncertainty of the stochastic gradient is higher with a smaller batch size, the results indicate that choosing a small $\overline{u}$ restricts potential

\begin{table}[h]\label{tab:results_different_u_epsilon2}
	\small
	\centering
	\caption{Numbers of iterations required for SGA-RMSProp to converge ($\varepsilon=\varepsilon_2$), reported as (mean, standard deviation) over three instances for $\overline{u}_i$, $\underline{u}_i$, $i=1,\ldots,4$.}
	\begin{tabular}{cccccc}
		\toprule
		Batch sizes & Problems & ($\underline{u}_1$, $\overline{u}_1$)& ($\underline{u}_2$, $\overline{u}_2$) & ($\underline{u}_3$, $\overline{u}_3$) &($\underline{u}_4$, $\overline{u}_4$) \\
		\midrule
		\multirow{12}{*}{$B=50$}&(ED, 20, 0.7) & (140, 1.0) & (177, 3.1) & (233, 7.9) & (288, 5.5)\\
		&(ED, 20, 0.2) & (144, 4.6) & (183, 4.4) & (253, 9.5) & (298, 5.9)\\
		&(ED, 50, 0.7) & (305, 11.9) & (364, 21.0) & (507, 26.7) & (601, 36.3)\\
		&(ED, 50, 0.2) & (309, 35.7) & (386, 36.0) & (504, 67.2) & (635, 49.7) \\
		&(ED, 100, 0.7) & (541, 2.3) & (634, 24.0) & (796, 31.8) & (1005, 13.1)\\
		&(ED, 100, 0.2) & (603, 18.2) & (765, 44.7) & (981, 56.9) & (1246, 43.7)\\
		&(AD, 20, 1) & (167, 14.2) & (191, 19.0) & (224, 28.6) &(297, 45.3) \\
		&(AD, 20, 2) & (154, 3.1) & (177, 7.5) & (215, 5.1) & (267, 15.9) \\
		&(AD, 50, 1) & (345, 49.6) & (378, 52.5) & (428, 72.9) &(512, 115.5) \\
		&(AD, 50, 2) & (331, 10.7) & (370, 16.4) & (427, 22.6) &(543, 37.6)  \\
		&(AD, 100, 1) & (568, 54.9) & (607, 53.9) & (650, 65.0) & (732, 73.9)\\
		&(AD, 100, 2) & (587, 44.3) & (631, 61.7) & (720, 74.7) &(866, 141.7)\\
		\midrule
		\multirow{12}{*}{$B=1000$}&(ED, 20, 0.7) & (143, 20.6) & (209, 29.6) & (281, 28.1) & (308, 8.7)\\
		&(ED, 20, 0.2) & (176, 19.7) & (246, 28.9) & (304, 8.3) & (324, 1.5)\\
		&(ED, 50, 0.7) & (299, 6.7) & (414, 10.1) & (603, 26.9) & (739, 29.5)\\
		&(ED, 50, 0.2) & (322, 66.8) & (482, 124.4) & (664, 126.6) & (784, 68.7) \\
		&(ED, 100, 0.7) & (485, 13.1) & (584, 41.8) & (777, 63.8) & (1093, 95.0)\\
		&(ED, 100, 0.2) & (537, 37.2) & (710, 128.5) & (1005, 258.0) & (1354, 321.7)\\
		&(AD, 20, 1) & (136, 31.5) & (81, 5.0) & (83, 9.0) & (89, 11.8) \\
		&(AD, 20, 2) & (116, 8.7) & (88, 2.9) & (92, 3.5) & (111, 3.6) \\
		&(AD, 50, 1) & (207, 21.1) & (180, 24.3) & (184, 23.8) & (193, 30.9)\\
		&(AD, 50, 2) & (245, 8.6) & (206, 3.1) & (215, 9.1) & (238, 5.5) \\
		&(AD, 100, 1) & (323, 29.0) & (293, 48.5) & (295, 50.0) & (304, 50.5)\\
		&(AD, 100, 2) & (366, 28.0) & (368, 29.8) & (384, 32.3) & (411, 38.9)\\
		\bottomrule
	\end{tabular}
\end{table}

\begin{table}[H]\label{tab:results_different_u_epsilon3}
	\small
	\centering
	\caption{Numbers of iterations required for SGA-RMSProp to converge ($\varepsilon=\varepsilon_3$), reported as (mean, standard deviation) over three instances for $\overline{u}_i$, $\underline{u}_i$, $i=1,\ldots,4$.}
	\begin{tabular}{cccccc}
		\toprule
		Batch sizes & Problems & ($\underline{u}_1$, $\overline{u}_1$)& ($\underline{u}_2$, $\overline{u}_2$) & ($\underline{u}_3$, $\overline{u}_3$) &($\underline{u}_4$, $\overline{u}_4$) \\
		\midrule
		\multirow{12}{*}{$B=50$}&(ED, 20, 0.7) & (198, 16.2) & (280, 6.6) & (336, 12.7) &(362, 1.7)\\
		&(ED, 20, 0.2) & (203, 13.3) & (287, 7.9) & (344, 5.9) & (357, 2.3)\\
		&(ED, 50, 0.7) & (380, 32.7) & (566, 63.9) & (738, 41.4) &(839, 5.5) \\
		&(ED, 50, 0.2) & (418, 69.8) & (567, 56.2) & (723, 58.4) &  (814, 18.2)\\
		&(ED, 100, 0.7) & (606, 9.8) & (890, 26.2) & (1196, 28.3) &(1503, 23.0) \\
		&(ED, 100, 0.2) & (753, 54.0) & (1054, 103.7) & (1320, 119.4) &(1540, 135.6) \\
		&(AD, 20, 1) & (194, 20.1) & (254, 33.6) & (341, 33.5) &(448, 67.5) \\
		&(AD, 20, 2) & (181, 15.0) & (237, 10.0) & (315, 11.6) & (395, 10.6) \\
		&(AD, 50, 1) & (377, 52.5) & (457, 88.0) & (552, 136.3) &(782, 205.6) \\
		&(AD, 50, 2) & (369, 22.5) & (458, 26.9) & (617, 20.7) &(846, 42.4) \\
		&(AD, 100, 1) & (606, 65.9) & (648, 85.8) & (771, 107.9) &(960, 143.4) \\
		&(AD, 100, 2) & (621, 51.8) & (697, 72.4) & (917, 107.3) &(1281, 207.6)\\
		\midrule
		\multirow{12}{*}{$B=1000$}&(ED, 20, 0.7) & (199, 34.0) & (294, 17.5) & (331, 9.5) & (335, 6.1)\\
		&(ED, 20, 0.2) & (241, 24.6) & (314, 18.6) & (343, 2.1) & (344, 3.1)\\
		&(ED, 50, 0.7) & (413, 23.1) & (575, 39.6) & (764, 31.1) & (807, 27.6) \\
		&(ED, 50, 0.2) & (431, 104.1) & (644, 152.5) & (807, 89.9) & (864, 21.2) \\
		&(ED, 100, 0.7) & (547, 50.5) & (720, 90.4) & (1110, 126.3) & (1484, 52.4)\\
		&(ED, 100, 0.2) & (666, 111.9) & (949, 280.3) & (1368, 300.7) &(1637, 146.3) \\
		&(AD, 20, 1) & (97, 3.1) & (85, 9.7) & (93, 10.0) & (114, 15.6)\\
		&(AD, 20, 2) & (108, 7.9) & (100, 3.1) & (120, 5.2) &  (155, 2.5)\\
		&(AD, 50, 1) & (200, 16.0) & (186, 23.1) & (201, 32.8) & (232, 60.4)\\
		&(AD, 50, 2) & (229, 11.9) & (218, 7.9) & (249, 20.0) & (328, 32.1)\\
		&(AD, 100, 1) & (324, 47.1) & (300, 42.2) & (314, 51.8) & (336, 62.8)\\
		&(AD, 100, 2) & (371, 29.8) & (388, 34.6) & (421, 35.2) &(480, 50.3)\\
		\bottomrule
	\end{tabular}
\end{table}
\noindent excessively large components in the descent direction $\boldsymbol{D}_k \boldsymbol{g}_k$ in \eqref{eq: vector form of x in RMSProp} that may result from the strong randomness in $\boldsymbol{g}_k$, thereby improving the stability of the algorithm.

On the other hand, the algorithm shows comparable performance across different combinations of $\overline{u}$ and $\underline{u}$ when $B=1000$. Moreover, the effect of $\varepsilon$ appears to be more significant than that of $\overline{u}$ and $\underline{u}$, so we focus on tuning $\varepsilon$, and set $\overline{u}$, $\underline{u}$ to $\overline{u}_2$, $\underline{u}_2$ respectively, which are not necessarily optimal for every experiment but are also not the worst choices.

Next, \cref{tab: results of different epsilon} concludes the numbers of iterations required for \hyperref[alg:SGA-RMSProp]{SGA-RMSProp} to converge with different $\varepsilon$ values and batch sizes. 
\begin{table}[h]\label{tab: results of different epsilon}
	\small
	\centering
	\caption{The numbers of iterations required for SGA-RMSProp to converge, reported as (mean, standard deviation) over three instances.}
	\begin{tabular}{@{}c@{\hspace{0.1cm}}c@{\hspace{0.1cm}}c@{\hspace{0.1cm}}c@{\hspace{0.2cm}}c@{\hspace{0.2cm}}c@{\hspace{0.1cm}}c@{\hspace{0.1cm}}c@{}}
		\toprule
		\multirow{2}{*}{Problems} & \multicolumn{3}{c}{$B=50$} & & \multicolumn{3}{c}{$B=1000$} \\
		\cmidrule{2-4} \cmidrule{6-8}
		 & $\varepsilon_1$ & $\varepsilon_2$ & $\varepsilon_3$ & & $\varepsilon_1$ & $\varepsilon_2$ & $\varepsilon_3$ \\
		\midrule
		(ED, 20, 0.7) & \textbf{(115, 0.6)} & (177, 3.1) & (280, 6.6) & & \textbf{(97, 3.1)} & (209, 29.6) & (294, 17.5) \\
		(ED, 20, 0.2) & \textbf{(113, 0.6)} & (183, 4.4) & (287, 7.9) & & \textbf{(100, 2.1)} & (246, 28.9) & (314, 18.6) \\
		(ED, 50, 0.7) & \textbf{(257, 3.2)} & (364, 21.0) & (566, 63.9) & & \textbf{(240, 3.1)} & (414, 10.1) & (575, 39.6) \\
		(ED, 50, 0.2) & \textbf{(266, 2.0)} & (386, 36.0) & (567, 56.2) & & \textbf{(244, 5.3)} & (482, 124.4) & (644, 152.5) \\
		(ED, 100, 0.7) & \textbf{(492, 2.0)} & (634, 24.0) & (890, 26.2) & & \textbf{(463, 3.5)} & (584, 41.8) & (720, 90.4) \\
		(ED, 100, 0.2) & \textbf{(550, 13.5)} & (765, 44.7) & (1054, 103.7) & & \textbf{(479, 13.1)} & (710, 128.5) & (949, 280.3) \\
		(AD, 20, 1) & \textbf{(153, 9.5)} & (191, 19.0) & (254, 33.6) & & (127, 5.1) & \textbf{(81, 5.0)} & (85, 9.7) \\
		(AD, 20, 2) & \textbf{(142, 0.6)} & (177, 7.5) & (237, 10.0) & & (110, 8.7) & \textbf{(88, 2.9)} & (100, 3.1) \\
		(AD, 50, 1) & \textbf{(333, 39.9)} & (378, 52.5) & (457, 88.0) & & (191, 18.4) & \textbf{(180, 24.3)} & (186, 23.1) \\
		(AD, 50, 2) & \textbf{(311, 6.4)} & (370, 16.4) & (458, 26.9) & & (206, 9.9) & \textbf{(206, 3.1)} & (218, 7.9) \\
		(AD, 100, 1) & \textbf{(559, 53.4)} & (607, 53.9) & (648, 85.8) & & (316, 16.9) & \textbf{(293, 48.5)} & (300, 42.2) \\
		(AD, 100, 2) & \textbf{(565, 45.5)} & (631, 61.7) & (697, 72.4) & & \textbf{(360, 24.7)} & (368, 29.8) & (388, 34.6) \\
		\bottomrule
	\end{tabular}
\end{table}

The results show that, with a small batch size of $B=50$, the algorithm performs better with $\varepsilon_3$ for all problems. Furthermore, when $B=1000$, the algorithm performs better with $\varepsilon_3$ on problems using ED, and performs better with $\varepsilon_2$ on problems using AD. These parameter settings will be used in future experiments.

In addition, the results from $B=1000$ indicate that a suitably larger $\varepsilon$ may accelerate the convergence of \hyperref[alg:SGA-RMSProp]{SGA-RMSProp}. Recalling \cref{thm: convergence of consistent LLSP}, the convergence rate $\gamma$ depends not only on $\left(1-\frac{2 \left(\underline{u} \lambda_d\left(\boldsymbol{A}^{\top}\boldsymbol{A}\right)-\overline{u} \sigma\right)}{ \overline{u}\left(\lambda_1\left(\boldsymbol{A}^{\top}\boldsymbol{A}\right)+\lambda_d\left(\boldsymbol{A}^{\top}\boldsymbol{A}\right)\right)}+\varepsilon\right)$ but also on $\text{G}(\gamma_1,\ldots,\gamma_K)$. Therefore, selecting an appropriate $\varepsilon$ may reduce $\text{G}(\gamma_1,\ldots,\gamma_K)$, resulting in faster convergence.

\subsubsection{Linear convergence rate}
\begin{figure}[H]\label{pic: SGA-RMSProp all convergence}
	\centering
	\includegraphics[width=0.9\textwidth]{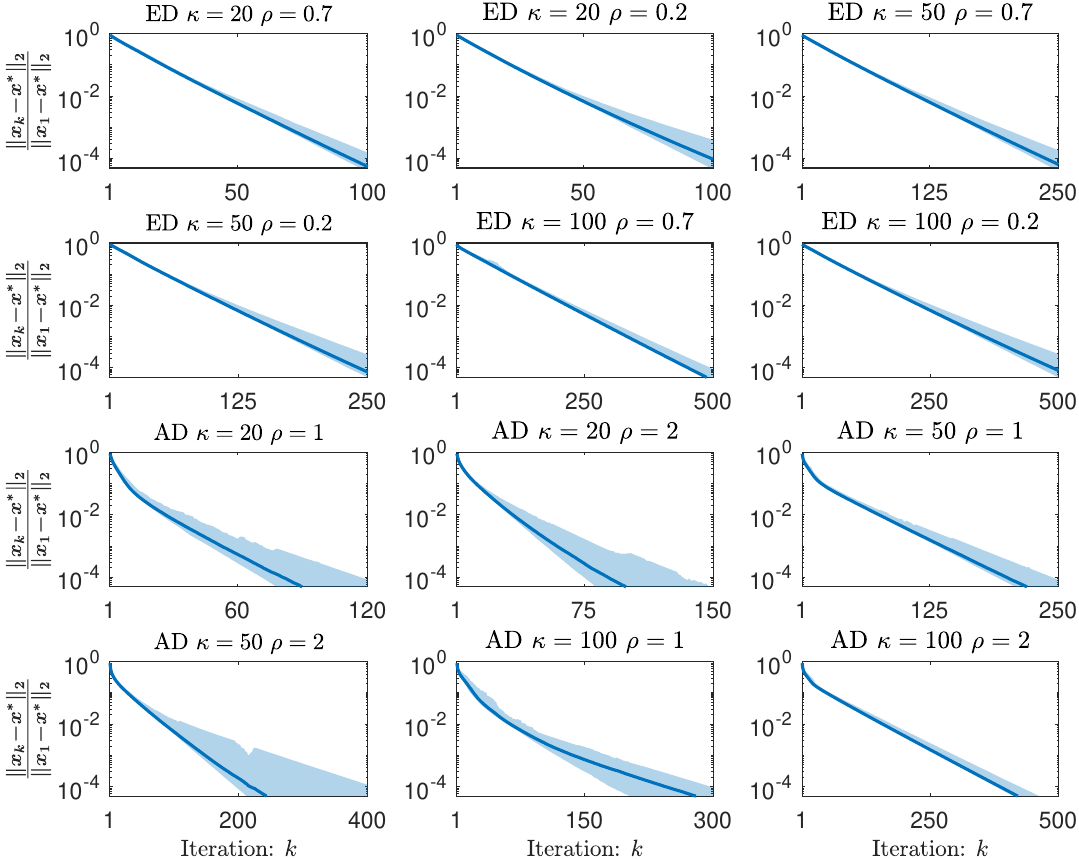}
	\caption{The R-linear convergence of SGA-RMSProp with $B=1000$.}
\end{figure}

We present the convergence of \hyperref[alg:SGA-RMSProp]{SGA-RMSProp} under different problems with $B=1000$ through \cref{pic: SGA-RMSProp all convergence}. The vertical axis represents $\frac{\|\boldsymbol{x}_k-\boldsymbol{x}^*\|_2}{\|\boldsymbol{x}_1-\boldsymbol{x}^*\|_2}$, while the horizontal
axis denotes the number of iterations. The dark line represents the average of the results in each iteration. The light shaded area shows the range from the 5th to 95th percentile of the results. Since the vertical axis is logarithmic scale and the plots eventually maintain a straight line under different problems, the results demonstrate that \hyperref[alg:SGA-RMSProp]{SGA-RMSProp} converges R-linearly on the consistent LLSP.

\vspace{-10mm}\subsubsection{Comparison of SGA-RMSProp and other algorithms}
We first compare our algorithm with the original RMSProp, whose iteration schemes follow \eqref{eq: RMSProp u iter}-\eqref{eq: RMSProp x iter} with all hyperparameters set as constants, on 12 instances derived from the problems (ED, 20, 0.7), (ED, 100, 0.7), (AD, 20, 2), and (AD, 100, 2). These instances cover both small and large condition numbers, as well as varying distributions of singular values \cite{strakovs1991real}. Since there is no widely accepted parameter choice for RMSProp on LLSP, we tune its parameters as best as we can. According to the documentation of PyTorch \cite{paszke2019pytorch}, the recommended default parameters for RMSProp are $\beta_k \equiv \beta = 0.99$ and $\eta_k\equiv \eta=0.01$. Based on these, \cref{tab:ori_RMSP_compare} reports the performance of RMSProp under combinations of $\beta=0.9,0.99,0.999$ and $\eta=0.01,0.1,1$, using the notation RMSProp($\beta$, $\eta$). In the experiments, we find that RMSProp is already slow when $\eta=0.01$, so smaller step sizes are not considered further. In \cref{tab:ori_RMSP_compare}, the column headings represent target accuracy level, defined as the relative error $\tfrac{\|\boldsymbol{x}_k-\boldsymbol{x}^*\|_2}{\|\boldsymbol{x}_1-\boldsymbol{x}^*\|_2}$ not exceeding $0.1$, $0.01$, $0.001$, and $0.0001$, respectively. Each cell contains a triple. The first component indicates the number of instances where the algorithm achieves the stated accuracy. The second component gives the mean number of iterations required to first hit the desired accuracy on those successful instances. The third component records the largest such iteration number among them. The triple (0, N, N) means that the algorithm does not reach the target on any instance.

\begin{table}[H]\label{tab:ori_RMSP_compare}
	\small
	\centering
	\caption{Comparison of RMSProp and SGA-RMSProp, reported as (number of successful instances, mean number of iterations, largest iteration number) over 12 instances.}
	\begin{tabular}{@{}c@{\hspace{2.5mm}}l@{\hspace{1mm}}r@{\hspace{2.5mm}}r@{\hspace{2.5mm}}r@{\hspace{2.5mm}}r@{}}
		\toprule
		Batch Sizes & Algorithm & $0.1$ & $0.01$ & $0.001$& $0.0001$ \\
		\midrule
		&RMSProp(0.9, 1) & (0, N, N) & (0, N, N) & (0, N, N) & (0, N, N)\\
		&RMSProp(0.9, 0.1) & (12, 104, 142) & (0, N, N) & (0, N, N) & (0, N, N)\\
		&RMSProp(0.9, 0.01) & (12, 692, 897) & (12, 843, 1134) & (0, N, N) & (0, N, N) \\
		&RMSProp(0.99, 1) & (6, 57, 64) & (0, N, N) & (0, N, N) &  (0, N, N)\\
		$B=50$&RMSProp(0.99, 0.1) & (12, 102, 172) & (6, 114, 147) & (5, 174, 209) &(3, 252, 260) \\
		&RMSProp(0.99, 0.01) & (12, 793, 1087) & (9, 929, 1109) & (6, 984, 1155) &(4, 979, 1153) \\
		&RMSProp(0.999, 1) & (12, 208, 386) & (9, 302, 621) & (9, 417, 840) &(7, 385, 1075) \\
		&RMSProp(0.9, 0.1) & (12, 100, 184) & (12, 241, 418) & (10, 341, 696) & (9, 411, 812) \\
		&RMSProp(0.999, 0.01) & (12, 528, 1107) & (8, 736, 1076) & (4, 777, 1149) & (3, 875, 925) \\
		&SGA-RMSProp & (12, 65, 121) & (12, 149, 257) & (12, 239, 431) &(12, 330, 607)\\
		\midrule
		&RMSProp(0.9, 1) & (0, N, N) & (0, N, N) & (0, N, N) &(0, N, N)\\
		&RMSProp(0.9, 0.1) & (12, 71, 122) & (0, N, N) & (0, N, N) & (0, N, N)\\
		&RMSProp(0.9, 0.01) & (12, 404, 486) & (12, 520, 598) & (0, N, N) &(0, N, N) \\
		&RMSProp(0.99, 1) & (2, 39, 40) & (0, N, N) & (0, N, N) &  (0, N, N)\\
		$B=1000$&RMSProp(0.99, 0.1) & (12, 63, 125) & (3, 70, 74) & (0, N, N) &(0, N, N) \\
		&RMSProp(0.99, 0.01) & (12, 509, 704) & (12, 743, 948) & (5, 774, 852) & (3, 897, 934) \\
		&RMSProp(0.999, 1) & (11, 127, 311) & (5, 172, 478) & (2, 155, 170) &(1, 187, 187) \\
		&RMSProp(0.999, 0.1) & (12, 71, 176) & (9, 189, 522) & (5, 234, 640) & (1, 176, 176) \\
		&RMSProp(0.999, 0.01) & (12, 256, 543) & (9, 490, 988) & (7, 625, 883) &(5, 749, 897) \\
		&SGA-RMSProp & (12, 48, 111) & (12, 113, 228) & (12, 183, 347) &(12, 256, 470)\\
		\bottomrule
	\end{tabular}
\end{table}

The results show that \hyperref[alg:SGA-RMSProp]{SGA-RMSProp} generally performs better than the original RMSProp. Therefore, in the following experiments, we do not include comparisons between SGA-RMSProp and RMSProp, but instead focus on comparing \hyperref[alg:SGA-RMSProp]{SGA-RMSProp} with SGD.

 Note that the SGD in our experiments also uses the mini-batch stochastic gradient \eqref{eq: k-th stochastic gradient}. We use the following suggested step size \cite{moorman2021randomized} for SGD when $B=50$:
\begin{equation*}
	\eta= \begin{cases}\frac{B}{\|\boldsymbol{A}\|_{\mathrm{F}}^2+(B-1) \lambda_d(\boldsymbol{A}^{\top}\boldsymbol{A})}, & \|\boldsymbol{A}\|_{\mathrm{F}}^2-(B-1)\left(\lambda_1(\boldsymbol{A}^{\top}\boldsymbol{A})-\lambda_d(\boldsymbol{A}^{\top}\boldsymbol{A})\right) \geq 0, \\ \frac{2B}{\|\boldsymbol{A}\|_{\mathrm{F}}^2+(B-1)\left(\lambda_1(\boldsymbol{A}^{\top}\boldsymbol{A})+\lambda_d(\boldsymbol{A}^{\top}\boldsymbol{A})\right)}, & \|\boldsymbol{A}\|_{\mathrm{F}}^2-(B-1)\left(\lambda_1(\boldsymbol{A}^{\top}\boldsymbol{A})-\lambda_d(\boldsymbol{A}^{\top}\boldsymbol{A})\right)<0.\end{cases}
\end{equation*}
Under our problems with the large batch size $B=1000$, we find that SGD generally performs better with $\eta=\frac{2}{\lambda_{1}(\boldsymbol{A}^{\top}\boldsymbol{A})+\lambda_{d}(\boldsymbol{A}^{\top}\boldsymbol{A})}$, the optimal fixed step size  of gradient descent for the strongly convex quadratic function \cite{Bertsekas99}. Therefore, this step size is used when $B=1000$. 

\begin{figure}[h]\label{pic:initial acceleration}
	\centering
	\includegraphics[width=0.9\textwidth]{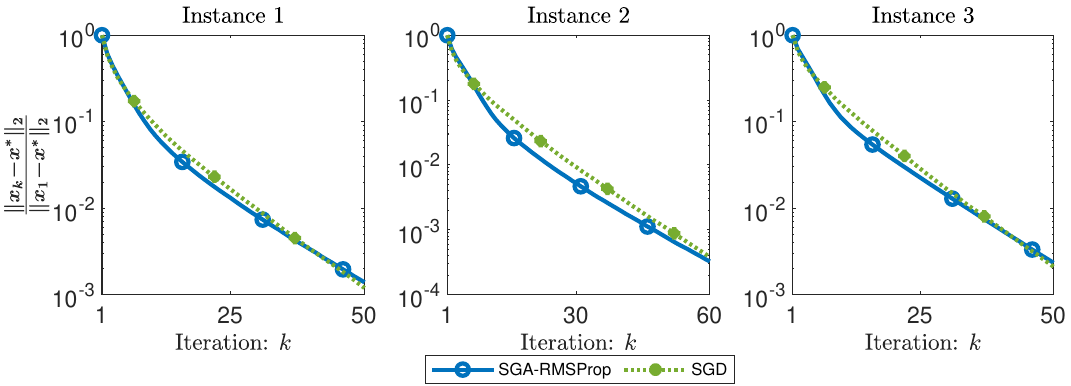}
	\caption{Comparison of SGA-RMSProp and SGD on (AD, 20, 1) with $B=1000$.}
\end{figure}

We directly compare the wall-clock time of the algorithms in \cref{tab: time of SGA-RMSProp and SGD}.
\begin{table}[h]\label{tab: time of SGA-RMSProp and SGD}
	\small
	\centering
	\caption{Wall-clock time (in seconds) of SGA-RMSProp and SGD, reported as (mean, standard deviation) over three instances.}
	\begin{tabular}{cccccc} 
		\toprule
		\multirow{2}{*}{Problems}  & \multicolumn{2}{c}{$B=50$} &  & \multicolumn{2}{c}{$B=1000$} \\
		\cmidrule{2-3} \cmidrule{5-6}
		  & SGA-RMSProp & SGD & & SGA-RMSprop & SGD \\
		\midrule
		(ED, 20, 0.7) & \textbf{(0.983, 0.002)} & (1.012, 0.005) & & (0.992, 0.027) & \textbf{(0.948, 0.003)}  \\
		(ED, 20, 0.2) & \textbf{(0.970, 0.006)} & (0.990, 0.004) & & (1.031, 0.013) & \textbf{(0.945, 0.005)}  \\
		(ED, 50, 0.7) & \textbf{(2.238, 0.028)} & (2.297, 0.040) & & (2.460, 0.053) & \textbf{(2.356, 0.014)}  \\
		(ED, 50, 0.2) & (2.266, 0.027) & \textbf{(2.193, 0.009)} & & (2.527, 0.049) & \textbf{(2.365, 0.013)}  \\
		(ED, 100, 0.7) & \textbf{(4.273, 0.002)} & (4.310, 0.014) & & (4.759, 0.086) & \textbf{(4.658, 0.031)}  \\
		(ED, 100, 0.2) & (4.689, 0.106) & \textbf{(4.210, 0.068)} & & (4.913, 0.134) & \textbf{(4.710, 0.004)}  \\
		(AD, 20, 1) & \textbf{(1.310, 0.068)} & (1.323, 0.082) & & (0.808, 0.064) & \textbf{(0.771, 0.033)}  \\
		(AD, 20, 2) & (1.223, 0.010) & \textbf{(1.219, 0.023)} & & (0.887, 0.013) & \textbf{(0.804, 0.011)}  \\
		(AD, 50, 1) & \textbf{(2.863, 0.377)} & (2.972, 0.380) & & \textbf{(1.825, 0.201)} & (1.865, 0.090)  \\
		(AD, 50, 2) & \textbf{(2.658, 0.048)} & (2.798, 0.052) & & (2.140, 0.040) & \textbf{(1.951, 0.041)}  \\
		(AD, 100, 1) & \textbf{(4.865, 0.499)} & (5.004, 0.657) & & \textbf{(2.992, 0.515)} & (3.131, 0.294)  \\
		(AD, 100, 2) & \textbf{(4.876, 0.376)} & (5.133, 0.430) & & (3.751, 0.291) & \textbf{(3.732, 0.239)}  \\
		\bottomrule
	\end{tabular}
\end{table}
The results indicate that \hyperref[alg:SGA-RMSProp]{SGA-RMSProp} performs better with the small batch size. Intuitively, when the number of samples is small, the stochastic gradient may not be accurate enough, whereas the correction step in \hyperref[alg:SGA-RMSProp]{SGA-RMSProp} helps the algorithm find a better descent direction. When the batch size is large, we observe that \hyperref[alg:SGA-RMSProp]{SGA-RMSProp} is faster than SGD in the early stages in terms of iterations, but is surpassed by SGD eventually. \cref{pic:initial acceleration} shows the early performances of \hyperref[alg:SGA-RMSProp]{SGA-RMSProp} and SGD on the three instances of the problem (AD, 20, 1) when $B=1000$.

Recalling the definition of $\gamma_k$ from the proof of \cref{thm: convergence of consistent LLSP}, $\beta_k<1$ may lead to $\gamma_k<1$, $k=1,\ldots,K$. This scenario frequently occurs during the early iterations of the algorithm, resulting in smaller values of $G(\gamma_1,\ldots,\gamma_k)$ at beginning, which may be the reason behind the fast initial convergence rate observed in our experiments. Therefore, starting with \hyperref[alg:SGA-RMSProp]{SGA-RMSProp} and switching to SGD at an appropriate time may improve the performance of the algorithm. Based on this idea, we conduct the next experiment.

\subsection{RMSP2SGD}
We introduce \hyperref[alg:RMSP2SGD]{RMSP2SGD} as follows, a method that enhances \hyperref[alg:SGA-RMSProp]{SGA-RMSProp} by adding an adaptive switching rule to SGD. The rule is that if $\beta_k$ equals 1 for five consecutive times, then the algorithm switches to SGD. As a result, it benefits from the initial rate of \hyperref[alg:SGA-RMSProp]{SGA-RMSProp} and SGD's faster convergence in the later stages.

We compare the wall-clock time of the three algorithms on the instances generated in \cref{subsubsec:Data description} with $B=1000$ and the parameters discussed in \cref{subsec:experiments for consistent LLSP}.

\cref{tab:time of three algorithms} indicates that \hyperref[alg:RMSP2SGD]{RMSP2SGD} improves the performance of \hyperref[alg:SGA-RMSProp]{SGA-RMSProp} on LLSP. It is also generally faster than SGD with respect to the wall-clock time, except in (ED, 20, 0.2). In addition, tuning the parameters may improve the initial convergence rate of \hyperref[alg:SGA-RMSProp]{SGA-RMSProp}, thereby making \hyperref[alg:RMSP2SGD]{RMSP2SGD} faster. This will be one of the future work.

\begin{algorithm}[H]\label{alg:RMSP2SGD}
	\caption{RMSP2SGD}
	\begin{algorithmic}[1]
		\REQUIRE Step size $\{\eta_k\}$, adjusted level $\varepsilon$, lower and upper bounds $\underline{u}$, $\overline{u}$, initial value $\boldsymbol{x}_1$.
		\STATE $\boldsymbol{u}_{0}=\left(\frac{1}{\overline{u}^2},\ldots,\frac{1}{\overline{u}^2}\right)^{\top}$.
		\FOR{$k=1,\ldots,K$}
		\STATE Sample a stochastic gradient $\boldsymbol{g}_k$.
		\STATE $\beta_k=\text{ \textbf{$\boldsymbol{\beta}$-Selection}}(\boldsymbol{g}_k,\boldsymbol{u_{k-1}},\varepsilon,\underline{u},\overline{u})$.
		\IF{$\beta_k$ equals 1 for five consecutive times}
		\STATE Break and switch to \textbf{SGD}.
		\ENDIF
		\STATE $\boldsymbol{u}_k=\beta_k \boldsymbol{u}_{k-1}+\left(1-\beta_k\right) \boldsymbol{g}_k^2$.
		\STATE $\boldsymbol{x}_{k+1}= \boldsymbol{x}_k-\frac{\eta_k}{\sqrt{\boldsymbol{u}_k}} \circ \boldsymbol{g}_k$.
		\ENDFOR
		\ENSURE The optimal point calculated by the algorithm $\boldsymbol{x}_{K+1}$.
	\end{algorithmic}
\end{algorithm}

\begin{table}[H]\label{tab:time of three algorithms}
	\small
	\centering
	\caption{Wall-clock time (in seconds) of the three algorithms with batch size $B=1000$, reported as (mean, standard deviation) over three instances.}
	\begin{tabular}{cccc} 
		\toprule
		Problems & RMSP2SGD & SGA-RMSProp & SGD\\
		\midrule
		(ED, 20, 0.7) & \textbf{(0.938, 0.003)} & (0.992, 0.027) &  (0.948, 0.003)   \\
		(ED, 20, 0.2) & (0.946, 0.001) & (1.031, 0.013) &  \textbf{(0.945, 0.005)}   \\
		(ED, 50, 0.7) & \textbf{(2.325, 0.018)} & (2.460, 0.053) &  (2.356, 0.014)   \\
		(ED, 50, 0.2) & \textbf{(2.343, 0.020)} & (2.527, 0.049) &  (2.365, 0.013)   \\
		(ED, 100, 0.7) & \textbf{(4.611, 0.029)} & (4.759, 0.086) &  (4.658, 0.031)   \\
		(ED, 100, 0.2) & \textbf{(4.666, 0.018)} & (4.913, 0.134) &  (4.710, 0.004)   \\
		(AD, 20, 1.0) & \textbf{(0.754, 0.121)} & (0.808, 0.064) &  (0.771, 0.033)   \\
		(AD, 20, 2.0) & \textbf{(0.768, 0.047)} & (0.887, 0.013) &  (0.804, 0.011)   \\
		(AD, 50, 1.0) & \textbf{(1.747, 0.170)} & (1.825, 0.201) &  (1.865, 0.090)   \\
		(AD, 50, 2.0) & \textbf{(1.948, 0.043)} & (2.140, 0.040) &  (1.951, 0.041)   \\
		(AD, 100, 1.0) & \textbf{(2.977, 0.357)} & (2.992, 0.515) &  (3.131, 0.294)   \\
		(AD, 100, 2.0) & \textbf{(3.621, 0.320)} & (3.751, 0.291) &  (3.732, 0.239)   \\
		\bottomrule
	\end{tabular}
\end{table}

\subsection{Experiments for inconsistent LLSP}
\begin{figure}[H]\label{pic: comparison on inconsistence LLSP}
	\centering
	\includegraphics[width=1\textwidth]{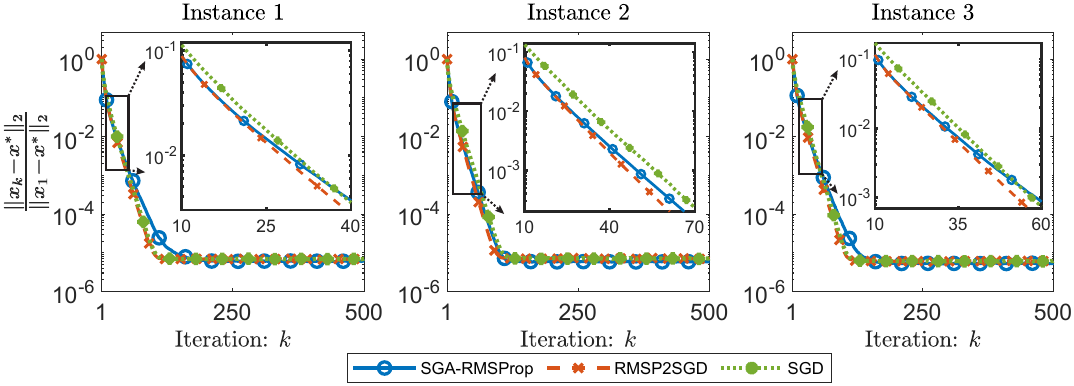}
	\vspace{-4mm}\caption{Comparison of the three algorithms for the inconsistent case.}
\end{figure}

In this subsection, we conduct experiments for the inconsistent LLSP. We generate the matrix $\boldsymbol{A}$ as described in \cref{subsubsec:Data description} and let $\tilde{\boldsymbol{x}}\in\mathbb{R}^d$ have independent standard normal entries, $\boldsymbol{b}=\boldsymbol{A} \tilde{\boldsymbol{x}}+\boldsymbol{\tau}$, where $\boldsymbol{\tau}$ represents a perturbation uniformly randomly drawn from a sphere with a radius of $10^{-3}$. We fix the problem at (AD, 20, 1) and present the results of the three instances. As shown in \cref{thm: convergence of inconsistent LLSP}, \hyperref[alg:SGA-RMSProp]{SGA-RMSProp} converges to a region of confusion in the inconsistent case, so we do not use the previous convergence criteria in this experiment. Instead, we set the maximum number of iterations to 500 and present the performances of the three algorithms in \cref{pic: comparison on inconsistence LLSP}.

 The results indicate that all three algorithms converge to the region of confusion. The inset plot of each figure shows that \hyperref[alg:SGA-RMSProp]{SGA-RMSProp} is still faster than SGD in the early stages, and \hyperref[alg:RMSP2SGD]{RMSP2SGD} is the fastest among three algorithms, which is similar to the consistent case.

\subsection{Real data}
\hspace{-1mm}We test the performance\hspace{-0.3mm} of \hyperref[alg:SGA-RMSProp]{SGA-RMSProp}\hspace{-0.3mm} on real data using the \hspace{-0.3mm}YearPredictionMSD dataset from UCI machine learning data repository \cite{misc_year_prediction_msd_203}. Previous studies used 

\vspace{-3mm}\noindent this dataset to predict song release years from timbre features through linear regression, and finally formulated it as LLSP \cite{bhardwaj2018ensmallen}. The dataset contains $n=463715$ records and $d=90$ timbre features. We apply mean-centered standardization to prepare the data for linear regression. Next, we set $\varepsilon$ to $\varepsilon_2$  and use a batch size $B=2000$, with all other parameters remaining the same as in the previous large batch size case.

\begin{figure}[H]\label{pic: realdata}
	\centering
	\includegraphics[width=0.5\textwidth]{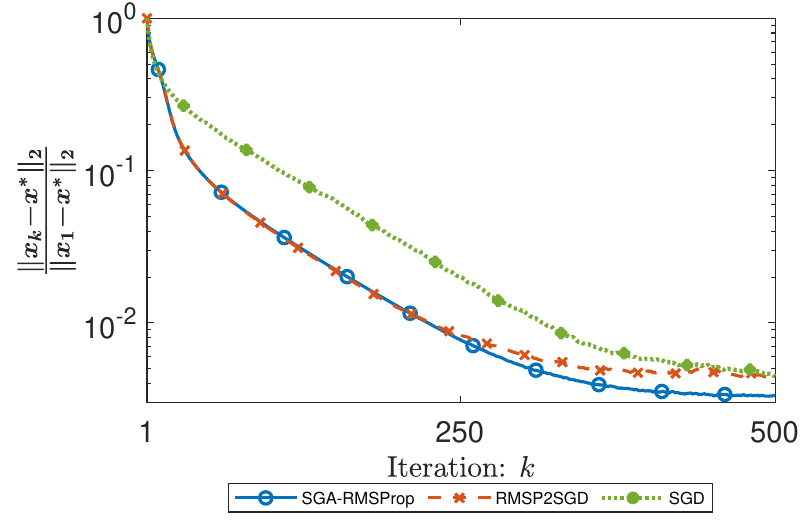}
	\vspace{-3mm}\caption{Comparison of the three algorithms on YearPredictionMSD.}
\end{figure}

\vspace{-2mm}The results are presented in \cref{pic: realdata}, indicating that the three algorithms converge to regions of confusion and \hyperref[alg:SGA-RMSProp]{SGA-RMSProp} achieves a smaller region. The plot also shows that \hyperref[alg:SGA-RMSProp]{SGA-RMSProp} is faster in the early stages. However, SGD can not surpass \hyperref[alg:SGA-RMSProp]{SGA-RMSProp} on this dataset, possibly due to entering the region of confusion too early.

\section{Concluding remarks}\label{sec:concluding remarks}
In this paper, we propose a method named stable gradient-adjusted RMSProp (\hyperref[alg:SGA-RMSProp]{SGA-RMSProp}), which uses a parameter called adjusted level to adaptively select the discounting factor and stably control the adjustment applied to the stochastic gradient. Theoretically, the R-linear convergence rate of \hyperref[alg:SGA-RMSProp]{SGA-RMSProp} is established on LLSP and a range of the adjusted level is provided to ensure this rate. Numerical experiments show that \hyperref[alg:SGA-RMSProp]{SGA-RMSProp} outperforms RMSProp, and achieves better performance than SGD when the batch size is small. Moreover, adaptively switching from \hyperref[alg:SGA-RMSProp]{SGA-RMSProp} to SGD generally improves the performance compared to using SGD alone.

Future studies might focus on analyzing the generally faster initial convergence rate of \hyperref[alg:SGA-RMSProp]{SGA-RMSProp} compared with SGD on LLSP from a theoretical perspective. For example, whether a tighter upper bound less than 1 for $\text{G}(\gamma_1,\ldots,\gamma_K)$ in \cref{eq: gamma} can be obtained? Other topics include exploring whether our method can be extended to general Adam-style algorithms or general functions.

\bibliographystyle{siamplain}
\bibliography{references}

\end{document}